\documentclass[11pt,a4paper,nofootinbib]{article}
\usepackage[dvips]{epsfig}
\usepackage{amsmath}
\usepackage{amssymb}
\usepackage[all]{xy}

\newcommand{\D}{\mathbb{D}}
\newcommand{\C}{\mathbb{C}}
\newcommand{\B}{\mathbb{B}}
\newcommand{\R}{\mathbb{R}}
\newcommand{\PP}{\mathbb{P}}
\newcommand{\T}{\mathbb{T}}
\newcommand{\Ss}{\mathbb{S}}

\newcommand{\Z}{\mathbb{Z}}
\newcommand{\Hh}{\mathbb{H}}
\newcommand{\lgra}{\longrightarrow}
\newcommand{\ra}{\rightarrow}
\newcommand{\dsps}{\displaystyle}
\newcommand{\smooth}{\mathcal{C}^{\infty}}
\newcommand{\exc}{\setminus}

\newcommand{\td}{\widetilde}

\def \pc{almost-complex}
\def \pce{almost-complex }
\def \ph{pseudo-holomorphic}
\def \phe{pseudo-holomorphic }
\def \hq{hyperbolic}
\def \hqe{hyperbolic }
\def \hle{holomorphic }
\def \hct{hyperbolicity}
\def \hcte{hyperbolicity }

\def \hl{holomorphic}
\def \st{structure}
\def \ste{structure }
\def \sp{symplectic}
\def \spe{symplectic }
\def \aK{almost-K{\"a}hler}
\def \aKe{almost-K{\"a}hler }

\def \Kae{K{\"a}hler }
\def \pp{\partial}
\def \f{\frac}
\def \dd{{\rm d}}

\def \mcl{\mathcal}
\def \be{\begin{equation}}
\def \ee{\end{equation}}
\def \beq*{\begin{equation*}}
\def \eeq*{\end{equation*}}
\def \ba{\begin{eqnarray}}
\def \ea{\end{eqnarray}}
\def \ba*{\begin{eqnarray*}}
\def \ea*{\end{eqnarray*}}
\def \dps{\displaystyle}

\def \area{{\rm area}}
\def \length{{\rm length}}
\def \proof{\underline{\it Dem.} : }
\def \endproof{\bigskip \bigskip $\blacksquare$}

\newtheorem{thm}{Theorem}[section]
\newtheorem{lemma}{Lemma}[section]
\newtheorem{prop}{Proposition}[section]
\newtheorem{defi}{Definition}[section]
\newtheorem{cor}{Corollary}[section]
\newtheorem{exam}{Example}[section]

\newtheorem{rmk}{Remark}[section]
\newtheorem{defi-thm}{Definition-Theorem}[section]

\title{\textbf{Floer homology, symplectic and complex hyperbolicities}}
\author{Anne-Laure \textsc{Biolley}}

\date{\today}

\begin{document}

\maketitle

\vspace{-.8cm}

\begin{center}
Department of Mathematics, University of Toronto,\\
Email: {\tt
alb\,@\,math.toronto.edu}

\end{center}

\medskip

\begin{abstract}
On one side, from the properties of Floer cohomology, invariant
associated to a symplectic manifold, we define and study a notion
of symplectic hyperbolicity and a symplectic capacity measuring
it. On the other side, the usual notions of complex hyperbolicity
can be straightforwardly generalized to the case of almost-complex
manifolds by using \phe curves. That's why we study
the links between these two notions of hyperbolicities when a
manifold is provided with some compatible symplectic and
almost-complex structures. We mainly explain how the non-symplectic
hyperbolicity implies the existence of \phe curves,
and so the non-complex hyperbolicity. Thanks to this analysis, we
could both better understand the Floer cohomology and get new
results on almost-complex hyperbolicity. We notably prove results
of stability for non-complex hyperbolicity under deformation of
the \pce structure among the set of the almost-complex
structures compatible with a fixed non-hyperbolic symplectic
structure, thus generalizing Bangert theorem  that
gave this same result in the special case of the standard torus.
\end{abstract}

\renewcommand{\thefootnote}{\fnsymbol{footnote}}
\setcounter{footnote}{0}
\footnotetext{
\phe curves - hyperbolicity.


Work supported by Ecole Polytechnique (France) and
University of Toronto (Canada).
}

\renewcommand{\thefootnote}{\arabic{footnote}}
\setcounter{footnote}{0}

Symplectic manifolds are naturally provided with
compatible almost-complex structure. Let us recall that a
symplectic structure and an almost-complex structure are
``compatible'' if they define a Riemannian metric on the manifold;
this metric is called {\bf almost-K\"ahler} (and is a
K\"ahler metric if the almost complex structure is integrable).
This naturally raises the issue of links between the \spe
properties and the \pce properties of the manifold when it is
provided with two compatible \st s. This analysis is likely to
provide a new and interesting approach to both these fields, and
in fact it has already been done. When an \pce \ste is given, you can
define pseudo-holomorphic curves which generalize the notion of
holomorphic curves in complex manifolds. Their introduction
by Gromov \cite{gromov} in \spe geometry led to an explosion of research in
this field, allowing to solve many
problems and to define new symplectic invariants such as the
Gromov invariant defined by
counting pesudoholomorphic curves.\\
This work fits in the framework of this analysis about links
between \spe and \pce properties, and more precisely, we tackle it
from the point of view of hyperbolicities notions.\\
From the complex point of view, there is currently great
interest around the various notions of complex hyperbolicity
(see \cite{kob}, \cite{lang}) -initially complex notions but that
can be straightforwardly generalized to the almost complex case-
which are not equivalent in general,
but which happen to be in the compact case. Their defintions are
all based on (pseudo)-holomorphic curves.\\
From the symplectic point of view, the properties of Floer
cohomology -symplectic invariant coming from Floer's new approach
to Morse theory, inspired by \phe curves-
has allowed me to define a notion of symplectic hyperbolicity.\\
We have studied the links between theses notions of \spe and \pce
hyperbolicities, using \phe curves as natural link. This will
certainly help and has already helped to better understand the
structure of Floer cohomology -which is at the heart of much
research today in symplectic geometry but also Hamiltonian
dynamics and mathematical physics- and to shed a new light on
complex hyperbolicity (through this approach of
\pce hyperbolicity that still remains  little explored).\\
\par
More precisely, {\bf Floer cohomology} comes from the generalization of
Morse theory to the case of a functional on an infinite dimensional space:
this
cohomology is generated by the $1$-periodic orbits of a Hamiltonian, which
are the critical points of a functional on the loop space associated to the
Hamiltonian. For closed manifolds, this has been studied in details and
it has
been proven that the Floer cohomology coincides with the usual cohomology.
Thanks to this result, the Arnold conjecture about the number of
$1$-periodic orbits of an Hamiltonian has been proven. A version of Floer
cohomology for compact manifold with contact type boundary has also been
defined, notably by Viterbo \cite{viterbo}, motivated by the Weinstein
conjecture about the existence of periodic orbits of an Hamiltonian $H$ on
a fixed energy level $\{H=a\}$, or equivalently about the existence
of Reeb orbits on a contact type hypersurface. Thus, in this
case,
the Floer cohomology is defined thanks to a sequence of Hamiltonians for
which the boundary of the manifold is an energy level (see section
\ref{floercohomreview} for more details). But, because of the result in the
closed case, one can naturally ask whether this Floer cohomology is
isomorphic to the usual one $H^*(M,\pp M)$.
This is the idea behind our definition of \spe
\hct. A compact manifold with contact type boundary (also called
$\omega$-convex) is said to be {\bf \sp-\hq} if the natural map between the
Floer cohomology and the usual one is surjective. We will explain this
in this paper. By construction, it will be noticed that if a manifold
is non-\sp-\hq,
then there exists a Reeb orbit on the boundary and the manifold
satisfies the Weinstein conjecture.\\
Moreover, we are going to define a {\bf \spe capacity} measuring \spe
\hcte (infinite if and only if the manifold is \spe \hq). And this will
allow us to naturally extend our notion of \spe \hcte to a larger class
of manifolds:
the open $\omega$-convex manifolds. An open \spe manifold is said
to be {\bf $\omega$-convex} if it can be written as an increasing union of
$\omega$-convex compact domains (\emph{i.e.} compact domains with contact
type boundary). These manifolds and their quotients will be the object
of study in this paper. We will often especially focus on the case of
{\bf Stein
manifolds} (in a much more general sense than the one of \cite{elgr}: see
definition \ref{steindefinition}), a particular class of $\omega$-convex
manifolds. For these $\omega$-convex manifolds, {\bf \spe \hct} will be
defined qualitatively as: the faster
the growth of the capacity of the compact domains, the more the manifold
is \sp-\hq. Depending on the quantity to which the growth of the capacity is
compared, different quantitative definitions of this notion can be written.\\
In this paper, after motivating and giving their definitions,
we study these notions of \spe \hcte and capacity, notably getting
some estimates on the growth of the capacity in the case of Stein
manifolds, and establishing some results on the \spe \hcte of product
manifolds.\\
Then, we will explain how
complex hyperbolicity implies \spe hyperbolicity, using as
a natural link \phe curves. \\
More precisely, the main results that will be presented here are:
 first a result linking non-\sp-\hcte and the existence of \phe disks:
\begin{thm}(theorem \ref{sympdisk})
Let $(M,\omega)$ be compact \spe manifold with contact type
boundary. If it is not \spe \hqe then, for any \pce \ste $J$
preserving the $\omega$-convexity (\emph{i.e.} preserving the
contact hyperplanes), there exists a $J$-\hle disk $f:\D \ra M$,
whose area is less than the capacity of $M$,
and whose boundary is included in the
boundary of $M$.
\end{thm}
Thus if an open $\omega$-convex manifold is non-\sp-\hq, this result provides
us with a sequence of bigger and bigger \phe disks with a control on their
area and the positions of their boundary, for any compatible \pce \ste
preserving the \spe convexity.\\
So secondly, we deduce that if a
$\omega$-convex manifold satisfies some properties of non-\sp-\hct, then
for any compatible \pce \ste $J$ preserving the $\omega$-convexity,
$(M,j)$ is not \pc-\hq, and more precisely it is not \aK-\hq, {\bf the
\aK-\hct} being one of the notion of complex \hcte that will be introduced
in section \ref{complexhyppart} and is equivalent with all the
others under some
compactness assumptions. This reads through different theorems, depending on
the non-\sp-\hcte assumption which is made, such as:
\begin{thm}
Let $(M^{2n},\omega)$ be an open $\omega$-convex \spe manifold
satisfying the non-\sp-\hcte assumption:
there exists an exhaustive increasing
sequence $(M_j)$ with $\mu(M_j)= O(d_g^2(\pp M_j, x_0))$. Then for any
(uniformly) compatible almost-complex \ste preserving the
$\omega$-convexity (of the $M_j$), $(M,J)$ is not \aK-\hq.
\end{thm}

\begin{thm}
Let $(M, J, \omega, \psi)$ be a Stein manifold and
$M_a=\{ \psi \leq a\}$. If $M$ is not \spe \hq, more precisely
if there exists an exhaustive sequence $(a_n)$ such that
$\mu(M_{a_n}) = O(a_n)$, then $M$ is not \aK-\hq.
\par
Moreover, if $M$ satisfies the stronger assumption of non-\sp-\hct:
$\mu(M_{a_n})= o(a_n)$, then $M$ is not
Kobayashi-\hq.
\end{thm}

By contrapositivity, if a manifold is complex-\hq, then we get a
lower bound
on the growth of the capacity of convex domains, \emph{i.e.} a result of
\sp-\hct, and thus we could better understand the Floer cohomology.
More precisely, we notably get:
\begin{thm}
Let $(M,\omega)$ be a $\omega$-convex \spe manifold. If there exists $J$
a compatible \pce \ste such that $(M,J)$ is \aK-\hq, then (if $g$ denotes
the \aKe metric associated with $\omega$ and $J$), $(M,\omega)$ is
symplectic \hq: for any increasing exhaustive sequence of $J$-convex
domains $(M_j)$,
\beq*
\f{\mu(M_j)}{d_g^2(\pp M_j, x_0)} \ra \infty.
\eeq*
\end{thm}

\begin{thm}
Let  $(M, J, \omega, \psi)$ a Stein manifold and  $M_a=\{ \psi
\leq a\}$.  If $M$ is \aK-\hq, then $(M,\omega)$ is  \spe \hqe:
\beq* \f{\mu(M_a)}{a} \ra \infty \eeq* \\
If $M$ is Kobayashi-\hqe then there exists a constant $C>0$ such
that \beq* \mu(M_{a}) \geq C a \eeq*
\end{thm}
As an application, this analysis provides us with some new examples of
non-complex-\hqe manifolds, notably among Stein manifolds (this notably
implies that the subcritical polarizations introduced in \cite{biran}
are not complex-\hq).\\
However, motivated by the contractibility of the set of compatible \pce
\st s, and keeping in mind the isolated result of Bangert, theorem
\ref{bangertthm} \cite{bangert} about
the non-\pc-\hcte of the standard torus, we wanted to improve this result
by getting rid of this assumption by asking the \pce \ste to preserve the
$\omega$-convexity. Thanks to some isoperimetric results, we prove that
for a manifold satisfying some hypothesis of bounded geometry
(see section \ref{maintheoremssection} for more details on these
assumptions), if it
satisfies some non-\sp-\hcte assumption, then for any
(uniformly) compatible \pce \st, the manifold is not complex-\hqe (again
more precisely in the \aKe meaning). This is described through several
theorems such as:
\begin{thm}
 Let $(M, \omega, \psi, J_0)$ be a Stein
manifold. If it satisfies one of the following properties:

\begin{enumerate}
\item  The geometry of $M$ is $II_\delta$-bounded. And
$M$ satisfies the non-\sp-\hcte assumption:
\begin{equation}
 \mu(M_{a})=O(a^{(1-\delta)\f{m}{2}}),
\end{equation}
with $m$ an isoperimetric value of $M$.

\item  The geometry of $M$ is $I_\delta$-bounded. And
$M$ satisfies the non-\sp-\hcte assumption:

\beq* \mu(M_{a})=O\left( a^{\min(1-\delta,\, (1+\delta)
\f{m}{2})}\right) \eeq*
\end{enumerate}
Then, for any $J$ uniformly compatible with $\omega$, $(M,
\omega,J)$ is not \aK-\hq. Thus if $(M, \omega)$ possesses a
compact quotient $(W, \omega_0)$, then for any $J$ compatible with
$\omega_0$, $(M, J)$ is not Brody-\hq.
\end{thm}

If the non-\sp-assumption satisfied by the manifold is less strong,
then we
couldn't prove the non-complex-\hcte for any compatible \pce \ste
anymore. However, we can prove the non-complex-\hcte for any
compatible \pce \ste in an open neighborhood of the standard one. This
can be expressed as:
\begin{thm}
 Let $(M, \omega, \psi, J_0)$ be a Stein
manifold. If the geometry of $M$ is $I_\delta$-bounded.
And if moreover $M$ satisfies the non-\sp-\hcte assumption:

\beq* \mu(M_{a})=O\left( a^{\min(1,\, (1+\delta) \f{m}{2})}\right)
\eeq*
Then, $(M, \omega, J_0)$ is not \aKe \hqe and this is also true on
a $\mathcal{C}^1$-neighborhood of $J_0$: for any compatible \pce
\ste $J$ satisfying  $|\dd \dd^c_J \psi|<C$ for a constant $C>0$\,
$(M, \omega,J)$ is not \aK-\hq.
\end{thm}
These results are particularly interesting since generally complex-\hcte
is an open property,  and so non-complex \hcte is an unstable property.
However, thanks to this result, we prove that (under some assumptions
of bounded geometry) if we restrict ourselves to the set
of \pce \st s compatible to a non-\sp-\hqe \st, then the non-complex
\hcte can be, in a certain meaning, stable: either all of these \pce
\st s are non-complex-\hq, or at least the set of non-complex-\hqe \st s
contains an open neighborhood around the standard one.
These results have numerous (and should have even more later) applications
in the study of \hct. For example this allows to prove the non-\hcte of
some Stein manifolds for non-standard \pce \st s. This applies notably to:
\begin{cor}
Let $W= \C \PP^ n \exc \{ k \textrm { hyperplans}\}$, with $k \leq n$,
be provided with its standard \spe and complex \ste $J_0$. Then for any
\pce \ste compatible $J$ in a $\mathcal{C}^1$-neighborhood of $J_0$,
$(W,J)$ is not complex \hq.
\end{cor}
This analysis applies also to the case of product manifolds providing us
with a whole class of Stein manifolds that not only are non-complex-\hqe
for any \pce \st, but moreover such that their product with any manifold
with bounded geometry satisfies this same property.
This implies for example:
\begin{cor}
Let $(N, \omega_N, J_N, \psi_N)$ be a Stein manifold with
$II_\delta$-bounded geometry and which admits an isoperimetric value
$m_N \geq 2$. Then the product manifold
$\C^n \times N$ is not \aK-\hqe for any \pce \ste $J$ uniformly
compatible with the \spe \ste $\omega=  \omega_0
\otimes \omega_N$ (for the Riemannian metric $g=g_0 \oplus g_N$).\\
Thus for any compact quotient $N_0$ of $(N,\omega_N,J_N)$, for any
\pce \ste $J$ compatible with $\omega_0 \otimes \omega_N$ on
$\T^{2n} \times N_0$, the \pce manifolds $(\T^{2n} \times N_0,J)$
and $(\C^n \times N, J)$ are not Brody-\hq.
\end{cor}
And in particular:
\begin{cor}
Let $N$ be a compact hyperbolic manifold, quotient of the \hqe ball
$(\D^n, \omega_{{\rm hyp}})$ (also denoted by $\Hh^n$) by a cocompact
group of \hle isometries. Then for any \pce \ste $J$ compatible with
the product \spe \ste $\omega_0 \oplus\omega_{{\rm hyp}}$ on
$\T^{2n} \times N$, the manifold $(\T^{2n} \times N,J)$ is not Brody-\hq:
there exists a non-constant $J$-\hle map $f:\C \ra \T^{2n} \times N$.
\end{cor}
But more general results are proved.\\
\par
Let's detail the plan of this paper. In a first part \ref{complexhyppart},
we recall the various notions of complex hyperbolicities
and its recent almost-complex approach and we introduce the notion of
\aK-\hcte which happens to be equivalent to all the other notions
of complex \hcte under compactness assumptions and will appear as a
natural link between complex and \spe \hct. Then \ref{floercohomreview}
is devoted to a brief review on Floer cohomology allowing to define
\spe \hcte and the capacity measuring it in \ref{symphyppart}. Then follows
a study of these notions with notably in \ref{hypproductpart},
the case of product or bundle manifolds.\\
After this introductive part to all these \hcte notions, we study the
links between
complex and \spe hyperbolicities, first making the link between \spe \hcte
and the existence of \spe disks in \ref{sympdiskpart}, then establishing
several \pce results, based on Nevanlinna theory, which allows to link
the existence of some \phe disks with non-complex \hct. These results
established in view of this present study could be useful in a more
general framework. Here, first, we use them to explain how non-\sp-\hcte
implies the
non-complex-\hcte of any compatible \pce \ste preserving the convexity in
\ref{maintheoremssection}. We present some immediate applications of
these results in \ref{applicationw}. Then we deepen the study to get rid
of this assumption of the \pce \ste preserving the convexity and thus to
get some stability results for non-complex-\hct. The main results are
enounced in \ref{enouncethebig} and some of their applications are given in
\ref{applicexsection}. However much more applications of these results
and this study are expected and will be presented in a later paper.
In order to prove them, in \ref{isopsection} we establish some
isoperimetric lemmas on \phe curves (and more generally on quasi-minimizing
currents). Finally, we use these results to prove the main theorems
and their applications.

\section{(Pseudo)-complex hyperbolicities}
\label{complexhyppart}
As it was explained in the introduction, the classical notions of
complex hyperbolicities (and their properties), which are based on
holomorphic curves (\cite{lang}, \cite{kob}) can be
straightforwardly generalized to the \pce case, by using \ph
curves. There are several non-equivalent definitions. The most
intuitive idea of the notion of complex hyperbolicites might be
given by the {\it Brody-hyperbolicity}: an \pce manifold is
Brody-hyperbolic if the only \phe curve fron the complex plane to
the manifold are the
constants.\\
The other main definition is: an \pce manifold is {\it Kobayashi-\hq} if
the Kobayashi pseudo-distance (defined through chains of \phe disks, see
\cite{lang} or \cite{kob} for more details) is a true distance.\\
In fact, it is proved (Brody theorem) that all the definitions of
(pseudo)-complex
hyperbolicities are equivalent in the compact case.\\
This \pce approach of complex hyperbolicities is quite new and one
of the first ones to have tackled this issue, especially from the
\spe point of view is Bangert \cite{bangert}: he studied the
hyperbolicities of the \pce torus and shows:
\begin{thm}[Bangert]
\label{bangertthm}
For any \pce \ste $J$ on the torus $\T^{2n}$ compatible with the standard
\spe \st,  $(\T^{2n}, J)$ is not complex \hq.
\end{thm}
(All the notions of complex hyperbolicities are equivalent in this case).
The torus provided with its standard complex \ste is obviously not
complex \hqe since there exists a \hle plane. This result tells us that
the torus is still not complex \hqe (and so that there still exists a \phe
plane) for any \pce \ste compatible with the standard \spe \st.
This is particularly relevant
since complex \hcte is an open property and thus non-complex-\hcte
is usually a non-stable property. But this Bangert's result claims that
if you restricts yourself to the set of \pce \st s compatible to the
standard \spe \st, then the non-complex \hcte is stable ! \\

Through a \spe approach of the issues of complex hyperbolicities,
and more precisely using \spe \hct, we are going to generalize this
result to a
more general type of \spe manifolds (namely the non-\spe \hqe ones).\\

While trying to link \spe hyperbolicities with complex hyperbolicities,
there is a notion of complex \hcte that appears as a natural bridge
between the two: a notion of {\it \aKe \hct}:
\begin{defi}
Let $(M,J)$ be an \pce manifold provided with a compatible \spe \ste
$\omega$. It is {\bf \aKe-\hq} if there exists
a constant $C>0$ such that for any \phe disk $f :\D \ra M$,
$|f'(0)| <C $ (the metric being taken for the \aKe metric $\omega(.,J.)$).
\end{defi}
Contrary to to the other notions of complex hyperbolicities, this
notion depends not only on the \pce \ste but also on the \aKe
metric fixed on the manifold, and thus on the symplectic \st.
However, this notion appears naturally in the proof of Brody theorem about
the equivalence of all the different notions of hyperbolicities in the
compact case. So it happens to be equivalent to all the other
notions of hyperbolicities under the usual compactness assumptions.
More precisely let's give this usual Brody theorem, as it can be naturally
read in terms of this \aKe \hct:
\begin{thm}[Brody]
Let $(M,J)$ be an \pce manifold provided with a compatible \spe \ste
$\omega$. If it is not \aKe \hqe then notably,
\begin{itemize}
\item if $M$ is relatively compact in an \pce manifold $W$, then $\bar{M}$
is not Brody-\hq.
\item if $M$ is compact, or has a co-compact group of
\hle pseudo-isometries, then $M$ is not Brody-\hqe and so is not
Kobayashi-\hqe either. In fact in this case, all the definitions of complex
hyperbolicities are equivalent.
\end{itemize}
\end{thm}

In fact, in \cite{bangert}, even if he doesn't use
the terminology of \aK-\hct, Bangert proves the theorem above as a
corollary of the proposition:
\begin{prop}
For any \pce \ste $J$ on $\C^n$, uniformly compatible with the standard
\spe \ste (for the standard euclidian metric), $(\C^n,J)$ is not
\aK-\hq.
\end{prop}
We are going to generalize this proposition to more
general $\omega$-convex (and more precisely Stein) manifolds that are
not \spe \hqe (see mainly part \ref{deepeningpart}).\\
But let's first define and study this notion of \spe \hct.

\section{Floer cohomology and symplectic hyperbolicity}
\subsection{Short review of Floer cohomology}
\label{floercohomreview}
Let $(M,\omega)$ (satisfying $[\omega] \pi_2(M)=0$) be a
symplectic manifold, with contact type boundary $\pp M= \Sigma$.
Let us recall that
$M$ being of contact type boundary means that there exists a
vector field $X$ in the neighborhood of $\Sigma$, transverse to
$\Sigma$, and such that $\mathcal{L}_X \omega=\omega$. Therefore,
the $1$-form $\alpha$, defined in the neighborhood of $\Sigma$ as
$\alpha=i_\eta \omega$, satisfies $d\alpha=\omega$. Thus
$\sigma=\alpha_{|\Sigma}$ is a contact form on $\Sigma$ and
$\xi=\{\sigma=0\}$ defines contact hyperplanes.\\
One defines on $\Sigma$ the Reeb field $R$ as:
\begin{equation}
\begin{cases} i_R (\omega_{|T\Sigma}) = i_R d\sigma =0 \\
\sigma(R)=1.
\end{cases}
\end{equation}
The closed orbits of $R$ are called the Reeb orbits or the closed
characteristics of $\Sigma$. We define their action by
$\mathcal{A}(\gamma)=\int_{\gamma} \gamma^* \sigma$. It is equal
to $T$, the period of the orbit $\gamma$. We define the spectrum
$\mathcal{S}(\Sigma)$ of $\Sigma$ as the set of the actions of all
the closed characteristics of $\Sigma$. \\
One of the goals of this version of Floer cohomology is to study
Weinstein's conjecture about the existence of Reeb orbits. In some
way, this happens to be equivalent to the study of Hamiltonian
orbits on a fixed energy level $\Sigma=\{H=a\}$. Indeed, this
last question does not depend on the Hamiltonian but only depends on
the geometry of the hypersurface $\Sigma$: the line $\R X_H$ is
characterized by $i_{X_H} (\omega_{|T\Sigma}) =0$
($i_{X_H} \omega =\dd H$ and $H$ is constant on $\Sigma$) and so
coincides with the Reeb line (or
characteristic line) of $\Sigma$.\\
\par
Let $J$ be any compatible \pce \st, compatible with $\omega$ and
preserving its $\omega$-convexity, \emph{i.e.} preserving a
contact hyperplanes field of $\Sigma$. This implies that $\Sigma$
is $J$-convex:
\begin{defi}
Let $\Sigma$ be an oriented hypersurface
included in a almost complex manifold $(M,J)$.
Let us define $\zeta = T\Sigma \cap J T\Sigma$.
It is locally fully defined (including its orientation)
on $T\Sigma$ by $\{ \tau =0 \}$.
We say that $\Sigma$ is $J$-convex, or pseudo-convex, if $d\tau(v, Jv)
>0$ for all $v \in \zeta \exc \{0\}$.
\end{defi}
The $J$-convex surfaces have a very interesting property~:
the $J$-holomorphic curves can not be interiorly tangent to them. \\
In the present case, $\{ \sigma=0 \} = J \{ \sigma=0 \} \subset
T\Sigma \cap J T\Sigma$.
Considering the dimensions, this implies that simply
$\{ \sigma=0 \} = T\Sigma \cap J T\Sigma$.
Thus, we can choose $\tau= \sigma$ for the convexity definition.
And since $d\sigma = \omega$ on $\{ \sigma=0
\}$, we find that $\Sigma$ is $J$-convex as expected.\\

Let's fix this structure $J$ and the preserved contact hyperplanes;
let's call them $\xi$ as above. Without changing $\xi$,
$\sigma$, $R$ or $J$, we can suppose that $J R$ and $X$ are
colinear on $\Sigma$ (just by replacing $X$ by $X + a X_0$ with a
some real
valued function).  Thus $J R =f \, X$ with $f$ some real valued function.\\
\par

Let $\phi_t$ be the flow associated to $\eta$ in the neighborhood $V$
of $\Sigma$. It satisfies $\phi_t^* \alpha={\rm{e}}^t \alpha$.
For $\epsilon_0$ small enough, we can then define
\begin{equation}
\psi~: \begin{cases} \Sigma \times (1-\epsilon_0,1] \lgra V \\
(x,z) \lgra \phi_{\ln(z)}(x). \end{cases}
\end{equation}
It follows that $\psi ^* \alpha=z\sigma $ et $\psi ^* \omega=d(z\sigma)$.\\
Thus, possibly restraining the neighborhood $(V,\omega)$,
we can identify it to
$(\Sigma \times (1-\epsilon_0,1], d(z\sigma))$.
And so we can construct an extension of $(M, \omega)$,
$(\td{M}, \td{\omega})$~:
$\td{M}=M \cup \Sigma \times [1, + \infty)$,
by extending $\omega$ to~:
\begin{equation*}
\td{\omega}= \begin{cases} \omega \, {\rm{sur}} \, M \\
d(z\sigma) \, {\rm{sur}} \, \Sigma \times [1, + \infty). \end{cases}
\end{equation*}
We extend the map $z$, naturally defined as the second coordinate on
$\Sigma \times [1-\epsilon_0, \infty)$ (modulo the identification of
$V$ to $\Sigma \times (1-\epsilon_0, 1]$),
with a map on  $M \exc V$ that takes values in $[0, 1-\epsilon_0]$.\\

On $\Sigma \times (1-\epsilon_0, + \infty)$, we notice that
$T\td{M}=T\Sigma \oplus \R \f{\pp}{\pp z}$ and that $\td{\omega}= dz
\wedge \sigma + d \sigma$. Therefore, $\td{\omega} (R, v)=0$ if $v \in T
\Sigma$ and $\td{\omega}(R, \f{\pp}{\pp z})=-1$. This implies $i_R
\td{\omega}=-dz$.\\
About the \pce \st, we extend  $J$ to $\td{J}$ on $\td{M}$ in the
standard way making it independent from $z$ on $\Sigma \times [1, \infty)$),
so that the $\Sigma \times \{ z \}$ hypersurfaces are pseudo-convex
for $z \geq 1$. \\

In the following, we will note, for $1-\epsilon_0<a<b$, $M_{[a,b]}=
 \{ (x,z) \in
\Sigma \times (1-\epsilon_0, \infty) \, \mid \, z
\in [a,b]\}$, $S_a=\{(x,z) \mid z=a\}$ and $M_a=M\exc V \cup  \{ (x,z) \in
\Sigma \times (1-\epsilon_0, \infty) \, \mid \, z \leq a \}$.\\

\par
The Floer cohomology is built as a limit of Floer cohomologies associated to
a sequence of adapted Hamiltonians (for which the boundary is an energy
level). Let's explain this construction. \\

To any given Hamiltonian $H :\R \times \td{M} \lgra \R$ is associated
an action $\mathcal{A}_H(\gamma)$ : if $\gamma~:\Ss^1 \lgra \td{M}$ is an
arbitrary contractible loop, we define its action by~:
\begin{equation}
\mathcal{A}_H(\gamma)= \int_{\D} \bar{\gamma}^* \td{\omega} \, -
\int_{\Ss^1} H(t, \gamma(t)) dt,
\end{equation}
where $\bar{\gamma}$ is an extension of $\gamma$ to $\D$
(the integral is uniquely defined since $[\omega] \pi_2(M)=0$). \\
First, given $\gamma~: \Ss^1 \lgra \td{M}$ and $Y$ a variation
of $\gamma$
(\emph{i.e.} vector field along $\gamma$), we have:
\begin{equation*}
\begin{aligned}
\dd\mathcal{A}_H(\gamma)(Y) & = \int_{\Ss^1} -\gamma^*(i_Y \td{\omega})
-\gamma^*(dH(Y)) \\
  & = \int_{\Ss^1} \td{\omega} \left( \f{\pp \gamma}{\pp t} - X_H,Y \right)
\end{aligned}
\end{equation*}
Thus, the critical points of $\mathcal{A}_H$ are exactly the orbits of
$X_H$.\\
\par
We will only consider adapted Hamiltonians
\emph{i.e.} satisfying on $\{z \geq 1 \}$,
$H(t,(x,z))=h(z)$, where $h$ is a convex increasing map on
$[1, +\infty)$.
Then on $\{ z \geq a \}$, we have $dH=h'(z) dz$ and $X_H=h'(z) X$.
Consequently, if $\gamma$ is a $1$-periodic orbit
of $X_H$ included in $\{ z \geq 1 \}$, then it is included in some
$\Sigma \times \{ z_0 \}$
(since $X$ is tangent to $\Sigma$) and is thus of the form
$\gamma(t)=(\gamma_0(t),z_0)$, with $\dot{\gamma_0}(t)=h'(z_0)\,
X$. This implies that
$\td{\gamma}(t)= \gamma_0 \left( \f{t}{h'(z_0)}\right)$ is
a closed characteristic of $\Sigma$.\\
Reciprocally, if $\td{\gamma}$ is
a closed characteristic of $\Sigma$ with action $T$, and assuming that
there exists $z_0$ such that $h'(z_0)= T$, then
$\gamma(t)=(\td{\gamma}(T t), z_0)$ is a $1$-periodic orbit
of $X_H$.\\
Moreover, as $\td{\omega}=d(z \sigma)$, we have
\begin{equation}
\mathcal{A}_H
(\gamma)=\int_{\Ss^1} z_0 \, \gamma_0^* \sigma - h(z_0) \, = z_0
h'(z_0)-h(z_0).
\end{equation}
Thus we have a one to one correspondence between
\begin{itemize}
\item on one side, the closed characteristics of $\Sigma$ with action $T$,
such that there exists $z_0 \leq 1$ with $h'(z_0)=T$,
\item on the other side, the $1$-periodic orbits of $X_H$. Their action is
$z_0 h'(z_0)-h(z_0)$ with $h'(z_0)=T$.
\end{itemize}
Now let us remember the construction of the Floer cohomology (see
\cite{viterbo}). From now on, we will choose adapted Hamiltonians such that
$h'(z)$ is constant, on some set $\{z \geq A \}$, equal to $\lambda \notin
\mathcal{S}(\Sigma)$.
This way, the $1$-periodic orbits of $X_H$ will all be included in a
compact set. Moreover,
\begin{equation*}
\dd\mathcal{A}_H(\gamma)(Y)=\int_{\Ss^1} g \left( J \f{\pp
  \gamma}{\pp t} - J X_H, Y \right),
\end{equation*}
where $g$ is the metric
canonically associated to $\omega$ and $J$ ($g(u,v)=\omega(u, Jv)$).
Therefore, the gradient trajectories of $\mathcal{A}_H$ correspond
to maps
$u~: \R \times \Ss^1 \lgra \td{M}$ such that~:
\begin{equation}
\f{\pp u}{\pp s}=-J \f{\pp u }{\pp t} - \nabla H, \, \forall \, (s,t)
\in \R \times \Ss^1,
\end{equation}
that is to say, if we define $\bar{\pp}_J u = \f{\pp u}{\pp
  s} + J \f{\pp u }{\pp t}$, $\bar{\pp}_J u= - \nabla H$.
These maps are called {\bf Floer trajectories} associated to the adapted
 Hamiltonian $H$. It is known that (see \cite{her}) that the Floer
trajectories can not be tangent to the interior of surfaces $\Sigma
\times \{ z \}$ with $z \geq 1 $. This comes from the $J$-convexity of
these hypersurfaces (this is proved in \cite{her}, \cite{her2} when
$JX=a \, R$ with $a$ constant but we can check that the demonstration
can be easily adapted to our case $JX =f \, R$ with $f$ a real valued
function see \cite{alb}).
\par
Now let us define the following complexes~:
\begin{equation}
C^k_a(H)=\bigoplus_{x \in \mathcal{T}(k,a,b)} \Z_2 . x
\end{equation}
where  $\dsps \mathcal{T}(k,a)= \left\{ x \textrm{
$1$-periodic trajectory of $X_H$, } \nu(x,H)=-i_{CZ,H}(x)=k \textrm{ and
} \mathcal{A}_H(x) \geq a\, \right\}$ (and $i_{CZ}$ is the Conley-Zehnder
index). \\
Given $x$ and $y$ two periodic orbits, we consider the set of Floer
trajectories linking $x$ and $y$:
\begin{eqnarray}
\mathcal{M}(y,x)&=&\mathcal{M}(H,J,y,x)\\
&=&\{ u:\R \times \Ss^1 \ra
\td{M}, \bar{\pp}u=  \nabla H, \lim_{s \ra - \infty} u=y, \lim_{s \ra
  + \infty} u =x \}\nonumber
\end{eqnarray}
Under some generic hypothesis (if the $1$-periodic orbits are
non-degenerate and for a generic choice of $H$ and $J$), one can
show that $\dps \dim \mathcal{M}(y,x)=i_{CZ}(x)-i_{CZ}(y)
=\nu(y,H)-\nu(x,H)$. As each trajectory can be reparametrized
translating the variable $s \in \R$, one can prefer to use
$\td{\mathcal{M}}(x,y)=\mathcal{M}(x,y)_{/ \R}$.  Then, since
$\nu(y,H)=\nu(x,H)+1$, we obtain $dim \td{\mathcal{M}}(x,y)=0$.
Moreover, we can check that, for a generic $H$, since the
trajectories all stay in a  compact set due to the $J$-convexity,
$\td{\mathcal{M}}(x,y)$ is compact;
it is therefore a finite set.\\
Furthermore, if $u:\R \times \Ss^1 \ra \td{M}$, we define its energy:
\begin{equation*}
E(u)=\f{1}{2} \int_{- \infty}^{+\infty} \int_{\Ss^1} \left| \f{\pp u}{\pp s}
  \right|^2 + \left| \f{\pp u}{\pp t} - X_H \right|^2.
\end{equation*}
If $\dsps \lim_{s \ra - \infty} u(s,.)=y, \lim_{s \ra
  + \infty} u(s,.) =x$,
\begin{equation*}
\mathcal{A}_H(y)-\mathcal{A}_H(x)  =\int_{-
  \infty}^{+\infty} -\f{\pp}{\pp s}(\mathcal{A}_H(u(s,.))
 =\int_{- \infty}^{+\infty} \int_{\Ss^1}- \td{\omega}\left( \f{\pp u}{\pp t}
  - X_H, \f{\pp u}{\pp s}\right)
\end{equation*}
and we point out that
\begin{equation*}
E(u)=\f{1}{2} \int_{- \infty}^{+\infty} \int_{\Ss^1} \left| \f{\pp u}{\pp s}
+J \f{\pp u}{\pp t} - \nabla H \right|^2 \dd s \dd t + \mathcal{A}_H(y)-
\mathcal{A}_H(x).
\end{equation*}
Thus the minimum of energy for the trajectories linking $y$ and $x$ are
reached by the Floer trajectories. And finally, if $u \in \mathcal{M}(y,x)$:
\begin{equation}
\label{bigaction}
\begin{aligned}
E(u) = \mathcal{A}_H(y)-\mathcal{A}_H(x)&  = &
\int_{- \infty}^{+\infty} \int_{\Ss^1} \left| \f{\pp u}{\pp s}
  \right|^2
 = \int_{- \infty}^{+\infty} \int_{\Ss^1}  \left| \f{\pp u}{\pp s} -
X_H \right|^2 \\
& =& \f{1}{2} \int_{- \infty}^{+\infty} \int_{\Ss^1} \left| \f{\pp u}{\pp s}
  \right|^2 + \left| \f{\pp u}{\pp s} - X_H \right|^2.
\end{aligned}
\end{equation}
Thus, if there exists a Floer trajectory linking $x$ and $y$, then
$\mathcal{A}_H(y) \geq \mathcal{A}_H(x)$.\\
Consequently, we can define if $x \in C^k_a(H)$,
\begin{equation*}
\dps
\partial x = \sum_{y \in C^{k+1}_a(H)} {\rm Card}(\td{\mathcal{M}}(y,x))\,
. \, y
\end{equation*}
Thus, if we define $\dsps C^k(H,a,b)=C^k_a(H) {/C^k_b(H)}$,
$\partial: C^k(H,a,b)\ra C^{k+1}(H,a,b)$ is well defined.
And it has been proven that $\pp \circ \pp =0$.
This allows us to define the Floer cohomology $FH^*(H,a,b)$ of $H$
as the cohomology groups of the $(C^*(H,a,b), \pp)$ complex.
Several remarks can be made (more details are given in \cite{alb}):
\begin{itemize}
\item The complex $C^k(H,a,b)$ being defined as the quotient $C^k_a(H)
{/C^k_b(H)}$, the Floer cohomology satisfies the long exact sequence:
\begin{equation}
\label{exlong1}
FH^*(H,b,c) \ra FH^*(H,a,c)\ra  FH^*(H,a,b) \overset{\pp^*}{\ra}
 FH^{*+1}(H,b,c) \ra
\end{equation}
\item if $H_0$ and $H_1$ are two Hamiltonians with
$H_0\leq H_1$, then there exists a morphism $\Phi_{H_1, H_0}:
FH^*(H_1,a,b) \ra FH^*(H_0,a,b)$.
\item We have done the previous construction for $H$ $\smooth$.
Nevertheless, we point out that (if neither $a$ nor $b$ are actions of
orbits of $H_0$ and $H_1$) in the case that $H_0$ and $H_1$ are
$\mathcal{C}^0$-close enough, then there is an isomorphism
$FH^*(H_0,a,b)\simeq FH^*(H_1,a,b)$. Then, if $K$
denotes a continuous Hamiltonian, we define:
\begin{equation*}
FH^*(K,a,b)= \lim_{ H \ra K,
                   H \smooth}  FH^*(H,a,b)
\end{equation*}
\end{itemize}
For fixed $\lambda$, let us look at the Hamiltonian $K_\lambda$ which
vanishes on $M$ and
takes the value $\lambda (z-1)$ on $\Sigma \times [1, \infty)$.
Its orbits are, on one side the constant orbits of $M$ with a null action,
and on the other side the orbits associated to the closed characteristics of
$\Sigma$ with an action $T\in [T_0, \lambda]$ . \\
Indeed, if we consider a Hamiltonian $K_{\lambda, \epsilon}$
obtained by smoothing $K_\lambda$ on $[1-\epsilon,
1+\epsilon]$, its orbits are, in addition to the constants on $M$, the
ones associated to the closed characteristics of $\Sigma$ (of period $T$)
with an
action $zT-h_{\lambda,\epsilon}(z)$, with $z \in [1-\epsilon,
1+\epsilon]$ satisfying $ h'_{\lambda,\epsilon}(z)=T$ (and thus $T\in
[T_0, \lambda]$).
This action goes to $T$ when $\epsilon$ goes to 0 (since $z$ goes to 1
and $h_{\lambda,\epsilon}(z)$ goes to 0).\\
For $\lambda >\lambda'$, we have $K_\lambda > K_{\lambda'}$ and so
we can define the
Floer cohomology of $M$ as a projective limit (see \cite{viterbo}):

\begin{defi-thm}
For $a<b$, we define the {\bf partial Floer cohomology} of $M$ by the
projective limit:
\begin{equation*}
FH^*(M,a,b)=\lim_{\lambda \ra \infty} FH^*(K_\lambda, a,b)
\end{equation*}
Then the {\bf Floer cohomology} of $M$ is defined as the projective limit:
\begin{equation}
\label{defhom}
FH^*(M)=\lim_{\mu \ra \infty} FH^*(M,-\delta,\mu), {\textrm{ with }} \delta
\textrm{ any number  } >0.
\end{equation}
This is a \spe invariant independent of the choice of the \pce
\ste $J$
and invariant under \spe deformation $(\omega_t)$.\\
The cohomology groups $FH^*(M,-\delta,\mu)$ are called the
{\bf truncated Floer cohomology} and will be simply denoted by
$FH^*(M,\mu)$.
\end{defi-thm}
Thus, roughly, the Floer cohomology is generated
\begin{itemize}
\item on one side by the constant orbits (critical points of a small Morse
fonction) inside of $M$ whose action is close to zero
\item on the other side by the closed characteristics of $\Sigma$
whose action is its period (so $>T_0$).
\end{itemize}
Since the constant orbits inside $M$ generate the Morse cohomology of the
manifold, if we consider the truncated Floer cohomology for a small
$\delta$ (which
means, if we keep only the orbits whose action is close to zero) then
we get (see \cite{alb} for the sketch of the proof):
\begin{prop}
For any $\delta >0$ sufficiently small:
\begin{equation*}
H^{*+n}(M,\partial M) \simeq FH^*(M,-\delta, \delta) \simeq
FH^*(K_\lambda,-\delta, \delta).
\end{equation*}
\end{prop}
We also notice that:
\begin{prop}
For $\lambda > \mu$,
\begin{equation}
\label{rmqisom}
FH^*(K_\mu,-\delta, \infty) \simeq FH^*(K_\mu,-\delta, \mu) \simeq
FH^*(K_\lambda,-\delta, \mu) \simeq FH^*(M,-\delta,\mu)
\end{equation}
and thus,
\be
\label{invlim}
FH^*(M) = \lim_{\lambda \ra \infty} FH^*(K_\lambda,-\delta,\infty),
\ee
\end{prop}

The Floer cohomology of balls has been computed (see \cite{fh1}).
Notably:
\begin{exam}
\label{exhomdisque}
\beq*
FH^{n}(\B^{2n}(r),a,b)=
\begin{cases} \Z_2 \textrm{ if } a<0 <b< \pi r^2\\
0 \textrm{ otherwise}
\end{cases}
\eeq*
\beq*
FH^{n+1}(\B^{2n}(r),a,b)=
\begin{cases} \Z_2 \textrm{ if } 0<a <\pi r^2<b \\
0 \textrm{ otherwise}
\end{cases}
\eeq* Thus, $FH^n(B^{2n}(r))= 0$ and more generally $FH^*(B^{2n}(r))= 0$.
\end{exam}
We can deduce from this the Floer cohomology of the hyperbolic ball
$\B^{2n}(1)=\D^n \subset \C^n$ provided with the \hqe \spe form
$\omega_{{\rm hyp}}= \f{\omega_0}{(1-|z|^2)^2}$.
\begin{lemma}
\label{homdiskhyp} For any $r<1$,
\beq*
FH^{n}(\D^{n}(r),a,b, \omega_{{\rm hyp}})=
\begin{cases} \Z_2 \textrm{ if } a<0 <b< \pi \f{r^2}{1-r^2}\\
0 \textrm{ otherwise}
\end{cases}
\eeq*
\beq*
FH^{n+1}(\D^n(r),a,b, \omega_{{\rm hyp}})=
\begin{cases} \Z_2 \textrm{ if } 0<a <\pi \f{r^2}{1-r^2}<b \\
0 \textrm{ otherwise}
\end{cases}
\eeq*
\end{lemma}
By writing (\ref{exlong1}) for $H=K_\lambda$ with $a=-\delta$, $b=\delta$
and $c=\mu$, and by making $\lambda \ra \infty$, we get:
\begin{equation}
\label{exlong0}
FH^*(M,\delta, \mu)\ra  FH^*(M, -\delta, \mu) \ra H^{*+n}(M,\partial M)
\overset{\pp^*}{\ra} FH^{*+1}(M,\delta, \mu) \ra
\end{equation}
By letting $\mu$ go to infinity, one gets:
\begin{equation*}
FH^*(M,\delta, \infty)\ra  FH^*(M) {\ra} H^{*+n}(M,\partial M)  \overset{\pp^*}{\ra} FH^{*+1}(M,\delta, \infty) \ra
\end{equation*}
Let's denote $c_M^*$ the map between $FH^*(M)$ and $H^{*+n}(M,\pp M)$. Intuitively,
this morphism $c_M^*: FH^*(M) \ra H^{*+n}(M, \pp M)$ ``only keeps
the orbits whose action is close to zero'' (the constant orbits inside $M$).
It measures the difference between $FH^*(M)$ and the purely topological
invariant $H^*(M,\pp M)$.\\
According to the results for the Floer cohomology of closed
manifolds it is natural to ask whether $c^*_M$ is an isomorphism
or not ? Which means, what is the \spe information (not
purely topological) carried by the Floer cohomology. First, we can
notice that if $c^*_M$ is not an isomorphism, then there exists at
least one Reeb orbit (the Floer complex is not reduced to the
constant orbits). That's why this property of non-isomorphism has
been called by Viterbo {\it Algebraic Weinstein property}
\cite{viterbo}. This is the idea behind my definition of symplectic
hyperbolicity:
\subsection{Symplectic \hcte and \spe capacity}
\label{symphyppart}
\begin{defi}
A symplectic manifold, compact with contact type boundary, is said
{\bf \spe \hq}
if $c_M^*: FH^*(M) \ra H^{*+n}(M,\pp M)$ is surjective.
\end{defi}
The study of Viterbo in \cite{viterbo} provides us with some first results
and examples. One of the main example of applications is provided by the
Stein manifolds, (following \cite{biran}) in a much more general sense than
the usual one by \cite{elgr}, since
it doesn't include either the integrability of the \pce \ste or the
completude of the manifold.
\begin{defi}
\label{steindefinition}
A {\bf Stein manifold} is a \aKe manifold $(M,\omega,J)$ such that there
exists a proper exhaustive pluri-sub-harmonic function $\psi$ with
$\dd \dd^c_J \psi = \omega$. \\
A Stein domain is a domain $\{ \psi
\leq a \}$ of a Stein manifold $(M,\omega,J, \psi)$. A Stein
manifold (or domain) is sub-critical one can choose for $\psi$ a
subcritical function (a function whose every critical point is
subcritical).
\end{defi}
The easiest and most fundamental example of Stein manifolds is the standard
$(\C^n, i, \omega_0, \f{|z|^2}{4})$. Another classical one is the
hyperbolic ball  $(\D^n, \omega_{{\rm hyp}}, J_0, \phi)$,
with $\phi =\f{1}{4} \ln \left(\f{1}{1-|z|^2}\right)$.\\
\par
According to \cite{viterbo}, the subcritical Stein domains are \spe \hq.
Moreover, thanks to the transfer morphism, it is proved
that for $W \subset M$, a $\omega$-convex domain in $M$, compact manifold
with contact type boundary, if $(W, \omega)$ is \spe \hq, then $(M,\omega)$
is \spe \hq.\\
\par
We'd like to extend this notion of \spe \hct, defined only for compact
manifold with contact type boundary, to a more general context. For this
purpose, we are going to use a notion of capacity measuring the \spe \hcte of
compact domains. Let us explain its construction.\\
For any $\mu>0$ one can consider the truncated cohomology $FH^*(M,\mu)$ and
there is the morphism $c_M^*$ factorizes through:
\xymatrix{
FH^*(M) \ar[r] \ar[dr]_{c_M^*} &
FH^*(M, \mu) \ar[d]^{c^*_\mu} \\
 & H^{*+n}(M,\pp M) }
(roughly by keeping first only the orbits with action less than $\mu$ and
then the ones with action close to zero). For $\mu>0$ very small,
$FH^*(M,\mu)\simeq H^{*+n}(M,\pp M)$ and $c_\mu^*=id$ is an isomorphism. Thus,
one can define for $(M,\omega)$ compact symplectic manifold with contact
type boundary:
\begin{equation*}
\mu(M)=inf\{ \mu \mid \, \, FH^{n}(M,\mu) \overset{c^n_\mu}{\ra}
H^{2n}(M,\partial M) {\textrm{is not surjective}} \} \, \in \bar{\R}_+
\end{equation*}
\begin{lemma}
This defines a {\bf \spe capacity}. Moreover, $(M,\omega)$ is \spe \hqe
if and only if $\mu(M)=\infty$.
\end{lemma}
(See \cite{alb} for the proof). Thus, intuitively, the bigger the capacity is,
the more the manifold is \spe \hq. \\
Thanks to the results on their Floer cohomology, one can compute
the capacity for the standard euclidian balls and for the hyperbolic balls.

\begin{exam} With $\omega_0$ the standard \spe \ste on $\C^n$,
$\mu(\B^{2n}(r), \omega_0)= \pi r^2$ for any $r>0$.\\
For $\omega_{{\rm hyp}}$ the \hqe \spe \ste on the \hqe ball $\D^n$,
$ \dsps \mu(\D^n(r), \omega_{{\rm hyp}}) = \pi \f{r^2}{1-r^2}$
for any $r<1$.
\end{exam}

Our generalization of the notion of \spe \hcte is based on this
capacity. The notion being previously defined only for compact
manifolds with contact type boundary, the natural context for its
extension is the $\omega$-convex manifolds which can be written as
the increasing union of such domains:
\begin{defi}
An open \spe manifold $(M, \omega)$ is {\bf $\omega$-convex} if it can
be written as $M=\cup M_j$ with $(M_j)$ an increasing sequence of
$\omega$-convex domains, \emph{i.e.} compact domains with contact
type boundary.
\end{defi}
The Stein manifolds are a particular case of $\omega$-convex manifolds
(the compact domains $M_a =\{ \psi \leq a \}$ being $\omega$-convex).\\

The \spe capacities $\mu(M_j)$ of these domains are well-defined and then
naturally we define qualitatively:
\begin{defi}
If $(M,\omega)$ is a $\omega$-convex manifold, the faster the capacity of
the compact domains $\mu(M_j)$ grows,
the more $M$ is {\bf \spe \hq}.
\end{defi}
This definition will be clarified quantitatively by specifying the
growth of the capacity. \\
This way, one can get several different versions/definitions depending to
the quantity to which the capacity is compared with.\\
In the particular case of Stein manifolds $(M, \omega,J , \psi)$, the
natural sets to consider for the $M_j$ are the $M_a =\{ \psi \leq a \}$.
Thus, we will naturally study the growth of the capacity $\mu(M_a)$
with respect to $a$.
We could for example consider the linear growth as a {\it limit-growth}.
Then an example of \sp-\hcte assumption could read:

\begin{exam}
A Stein manifold is (Stein)-\spe \hqe (or $S$-\spe \hqe) if
$\dsps \lim_{a \ra \infty} \f{\mu(M_a)}{a} \ra \infty$.
\end{exam}
In the more general case of any $\omega$-convex \spe manifold
$M=\cup M_j$, one {\it limit-growth} that one could consider is the linear
growth in $d_g^2(x_0, \pp M_j)$, with $d_g$ the distance associated
to a metric $g$ (fixed compatible with $\omega$: $|\omega|_g=1$). Let's
notice that it is natural to compare the capacity to the
square of the distance since this
is a $2$-dimensional invariant. Thus the \sp-\hcte assumption could
for example read:
\begin{exam}
\label{examghyp} A $\omega$-convex manifold (provided with a
compatible Riemannian metric $g$) is ($g$-)\spe \hqe if for any
exhaustive increasing sequence of $\omega$-convex domains $(M_j)$,
$\dsps \lim_{j \ra \infty} \f{\mu(M_j)}{d_g^2(x_0, \pp M_j)} \ra
\infty$.
\end{exam}
In the following it will be convenient to name {\bf strongly \spe
\hq} the manifolds for which there exists a compact domain with
infinite capacity (then all the bigger domains have infinite
capacity too). In this case, the sequence of capacities of any
exhaustive sequence of domains will be trivially stationary equal
to infinity.\\
\par
Let us look at our basic examples. For the Stein manifold$(\C^n,
\omega_0, J_0, \f{|z|^2}{4})$, since $\mu(\B^{2n}(r))=\pi r^2$,
$\mu(M_a)= 4 \pi a$. Thus the growth of the capacity is linear and
with the definition above
$(\C^{2n}, \omega_0)$ is not $S$-\spe \hqe,\\
On the other side for the \hqe ball $(\D^n, \omega_{{\rm hyp}}, J_0,
\psi)$ with $\psi=  \f{1}{4} \ln\left(\f{1}{1-|z|^2}\right)$,
the sets $\{\psi \leq a \}$ are $\B^{2n}(1-e^{-4a})$ and so
$\mu(\{\psi \leq a \}, \omega_{{\rm hyp}}= \pi (1-e^{-4a}) e^{4a}$. Thus,
the growth of the capacity is exponential in $a$ and $(\D^n,
\omega_{{\rm hyp}})$ is $S$-\spe \hq.\\
\par
More generally, it is interesting to study the growth of the capacity of the
sets $M_a$ of Stein manifolds. \\
For this we can use their very special property: the gradient flow
of $\psi$ is conformal and satisfies $(\phi_t)^* \omega= e^t
\omega$. Moreover we check that if $\frac{||d\psi||^2}{\psi} \geq
\alpha$ on $\{ a \leq \psi \leq b \}$ then for $t=\frac{1}{\alpha}
log(\frac{b}{a})$, $\phi_{-t}(M_b) \subset M_{a}$. Thus, using
the conformality of the capacity, we get an
estimate of the growth of the capacity in terms of
$\frac{||d\psi||^2}{\psi}$ (see my thesis \cite{alb} for more details). More
precisely,
\begin{prop}
\label{resmu} If $(M,J, \omega, \psi)$ is a Stein manifold,
complete at infinity, and if there exists $\alpha>0$ and a compact
$K$ such that $\frac{||d\psi||^2}{\psi} \geq \alpha$ on $W \exc
K$, then, there exist constants $a_0>0$ $c >0$ such that for all
$a>a_0$
\begin{equation*}
\mu(M_a, \omega) \leq c a^{\f{1}{\alpha}}.
\end{equation*}
\end{prop}
This estimate can be refined in terms of two functions linked
with $\frac{||d\psi||^2}{\psi}$: let's define
\begin{defi}
\label{fctfond} Let $(M,J,\psi)$ be a Stein manifold ($\psi\geq
0$). For $s>0$~:
\begin{eqnarray*}
\alpha(s)= \inf_{ \{ \sqrt{\psi} =s \} } \f{|d \psi|^2}{\psi}
&\textrm{
  and }& \alpha_0(s)= \inf_{ \{ \sqrt{\psi} =s \} } \f{|d \psi
  |}{\sqrt{\psi}} \\
\beta(s)= \sup_{ \{ \sqrt{\psi} \leq s \} } \f{|d \psi|^2}{\psi}
&\textrm{
  and }& \beta_0(s)= \sup_{ \{ \sqrt{\psi} \leq s \} }
  \f{|d \psi |}{\sqrt{\psi}}
\end{eqnarray*}
\end{defi}
As explained in my thesis \cite{alb} these fonctionss are naturally linked
with the growth function introduced by Polterovich in \cite{polt}. In order
to better understand them and the linked results, let's have a look at
$\f{|d \psi|^2}{\psi}$ for the two fundamental examples of Stein
manifolds:\\
\begin{itemize}
\item For $(\C^n, \omega_0, J_0, \f{|z|^2}{4})$, $|\dd \psi|=\f{|z|}{2}$ and
$\f{|\dd \psi|^2}{\psi}=1$
\item For $(\D^n, \omega_{{\rm hyp}}, J_0,
\psi)$ with $\psi=  \f{1}{4} \ln\left(\f{1}{1-|z|^2}\right)$,
$|\dd \psi|=\f{|z|}{2}$ and $\f{|\dd \psi|^2}{\psi}=
\f{|z|^2}{\ln\left(\f{1}{1-|z|^2}\right)} \ra 0$ when $\psi \ra \infty$
(\emph{i.e.} $|z|\ra 1$).
\end{itemize}

Then using the same method as for proposition \ref{resmu}, we get:
\begin{prop}
Let $(M,\omega,J, \psi)$ be a Stein manifold with a finite number
of critical values, $M_a= \{ \psi \leq a \}$, and $\alpha$ and
$\beta$ the functions defined in \ref{fctfond}. Then there exists
a constant $c$ such that
\begin{equation}
\label{dansunsens} \mu(M_a) \leq c \,
{\rm{exp}}\left(\int_{a_0}^{a} \f{1}{\alpha(s)} \f{ds}{s} \right).
\end{equation}
In the same way,
\begin{equation}
\label{danslautre} \mu(M_a) \geq c \,
{\rm{exp}}\left(\int_{a_0}^{a} \f{1}{\beta(s)} \f{ds}{s} \right).
\end{equation}
\end{prop}

Thanks to these propositions we can study the \spe \hcte of Stein
manifolds: either for any (or equivalently for one) $a_0$
$\mu(M_{a_0})=\infty$ and $M$ is strongly \spe \hqe ; or there exists $a_0$
such that $\mu(M_{a_0})<\infty$ and the propositions above give us
a quite precise estimate of the growth of the capacity and so a
"measure of the \spe \hct" of the manifold.\\
For example, this provides us with a lower bound on the growth of the
capacity of Stein manifolds of the type $M=W \exc H$ with $H$ the zeros
of a pseudo-transverse section, see lemmas \ref{cappseudotrans1} and
\ref{cappseudotrans2}.

\subsection{Symplectic \hcte of product manifolds}
\label{hypproductpart}
In order to complete this study of \spe \hct, there is another
issue interesting to tackle: the \spe \hcte of product manifolds
(and fibered manifolds). Thanks to the results of A. Oancea
\cite{alex} on the Floer cohomology of product manifolds, we
prove:
\begin{prop}
Let $(M, \omega_M)$ and $(N, \omega_N)$ be two compact \spe
manifolds with contact type boundary. Then \beq* \mu(M\times N,
\omega_M \otimes \omega_N) \leq \min( \mu(M, \omega_M), \mu(N,
\omega_N)). \eeq*
\end{prop}
\proof According to \cite{alex}, there exists a spectral sequence
linking the Floer cohomology of $M^{2m} \times N^{2n}$ with the
one of $M$ and the one of $N$. More precisely, for any $\lambda
>0$, we have the commutative diagram:\\
 \xymatrix{ \bigoplus_{r+s =k} FH^r(M,-\delta, \lambda)
\otimes FH^s(N,-\delta, \lambda) \ar[r] \ar[d]_{c_M^r \times
c_N^s} & FH^k(M \! \times \! N,-\delta, 2 \lambda) \ar[d]_{c_{M \!
\times \! N}^k} \ar[r] &
FH^k(M \! \times \! N,-\delta, \lambda) \ar[d]_{c_{M \times N}^k}\\
\bigoplus_{r+s =k} H^{r+m}(M ,\pp M) \otimes H^{s+n} (N, \pp N)
\ar[r] & H^{k+m+n}(M \! \times \! N,\pp (M \! \times \! N))
\ar[r]^{\sim} & H^{k+m+n}(M \! \times \! N,\pp (M \! \times \! N))
}

 Moreover, roughly, the sequence of Hamiltonians needed to build the
 Floer cohomology of $M \times N$ "comes from" the sum of the
 Hamilonians for $M$ and the ones for $N$ (see \cite{alex} or
 \cite{alb} for
 more details) ; thus a periodic orbit for $M \times N$ with action
 less than $\lambda$ can only come from a periodic orbit of $M$
 and one of $N$ both with action less than $\lambda$. Thus,
 one can check that the composed above morphism:
 $$ \bigoplus _{r+s=k} FH^r(M,- \delta, \lambda) \otimes
FH^s(N,-\delta,\lambda)  \longrightarrow FH^k(M\times N, -\delta,
2\lambda) \longrightarrow FH^k(M\times N, -\delta, \lambda) $$ is
surjective. Furthermore, for maximal $k$, which means for $k=m+n$,
\begin{eqnarray*}
H^{k+m+n}(M\times N,\pp (M\times N)) & = & H^{2(n+m)} (M\times
N,\pp (M\times N))
= \Z_2 \\
\textrm{et } \bigoplus_{r+s =k} H^{r+m}(M ,\pp M) \otimes H^{s+n}
(N, \pp N) &=& H^{2m}(M, \pp M) \otimes H^{2n}(N, \pp N) = \Z_2
\end{eqnarray*}
are isomorphic. So at the end, if $FH^m(M, -\delta, \lambda)
\overset{c_M^m}{\ra} H^{2m}(M, \pp M)$ or
 $FH^n(N, -\delta, \lambda) \overset{c_M^n}{\ra} H^{2n}(N, \pp N)$
 is not surjective then $c^m_M \times c^n_N$ is neither. And so,
 according to the commutative diagram above,
$c^{m+n}_{M \times N}: FH^{m+n}(M\times N,-\delta, \lambda) \ra
H^{2(n+m)} (M\times N,\pp (M\times N))$ is not surjective.
\endproof \\
As a corollary,
\begin{cor}
The non-\sp-\hcte is stable under product with whatever other
$\omega$-convex \spe manifold.
\end{cor}
Let's explained this for example for the version given in example
\ref{examghyp} of \spe \hcte:\\
 \proof If $(M,\omega)$ is not ($g$-)\spe \hqe then there exists a
 sequence $(M_j)$ of $\omega$-convex domains of $M$ such that
 $\mu(M_j, \omega) \leq C d_g^2(x_0, \pp M_j)$. Let $(N, \omega')$
 be another $\omega$-convex manifold (provided with any metric
 $g'$). It is written as the union of $\omega$-convex domains
 $(N_j)$ and, possibly considering an extracted subsequence, one
 can assume that $d_{g'}(y_0, \pp N_j) \geq d_g(x_0, \pp M_j)$ for
 any $j$. Then if $\omega_0=\omega \otimes \omega'$ and $g_0$
 denotes the product metric on $M \times N$,
$\mu(M_j \times N_j, \omega_0) \leq \mu(M_j, \omega) \leq
C d_g(x_0,\pp M_j)^2 \leq C d_{g_0}((x_0,y_0), \pp(M_j \times
N_j))^2$.
\endproof  \\
Let s notice that this property is analogous to the one we get in
the framework of complex hyperbolicities: if $(M,J)$ is not complex
\hqe then for any almost-complex manifolds $(N,J')$, the manifold
$(M \times N, J \oplus J')$ is not complex \hqe either.\\
\par
We can adapt the same kind of reasoning for the fibered manifolds.
However, the Floer cohomology of fibered manifolds is still not
very well understood. We only have some precise results for some
\spe fiber bundles whose basis is closed and more precisely for
the line fibered bundle with negative curvature \cite{alex}. This
allows us to get:
\begin{prop}
Let $(B^{2n} \beta)$ be a closed \spe manifold with entire \spe
form: $[\beta] \in H^2(B, \Z)$ and $\mathcal{L}$ be line bundle
over $M$ such that $c_1(\mathcal{L}) = -[\beta]$. If $\mathcal{L}$
is provided with a {\it strong-fibered-\spe \ste} (\cite{alex})
$\omega$ then the capacity of $L_r$, the disk bundle included in
$\mathcal{L}$ over $B$ with radius $r$, satisfies $\mu(L_r,
\omega) \leq \mu(B(r), \omega_0) = \pi r^2$.
\end{prop}
\proof The spectral sequence that A. Oancea gets in \cite{alex}
and the morphism from it to the usual Leray-Spectral sequence
writes in maximal dimension (see \cite{alb} for more details):

\xymatrix{FH^{n+1}(L_r, -\delta, \mu) \ar[r]^{c^{n+1}_{L_r}}
\ar[d] &
H^{2(n+1)}(L_r,\pp L_r) \ar[d] \\
H^{2n}(B,FH^1(F,\mu)) \ar[r]_{c_F^1!} & H^{2n}(B,H^2(F, \pp F))}

with $c_F^1!$ being the morphism from $H^{2n}(B,FH^1(F,\mu))$ to
$H^{2n}(B,H^2(F, \pp F))$ naturally associated to the morphism:
$c_F^1: FH^1(F,\mu) \ra H^2(F, \pp F)$. Since $F \simeq
(\D^{2n}(r), \omega_0)$, for $\mu > \pi r^2$ $c_F^1$ is not
surjective, and even: $FH^1(F,\mu)=0$. Moreover $H^2(F \pp F) =
\Z_2$ and $H^{2(n+1)}(L_r,\pp L_r)) \ra H^{2n}(B,H^2(F, \pp F))$
is an isomorphism. Thus if $c_F^1!$ is not surjective then
$c^{n+1}_{L_r}$ is not either.
\endproof  \\
Now that both the context of complex and \spe hyperbolicities have
been settled and studied, let's study the links between the two.\\

As it was explained in the introduction, the \phe curves are used
as link between these two notions. They are naturally (by
definition) linked with the notion of (almost-)complex \hcte.
Let's study the link between \spe \hcte and the existence of \phe
curves.

\section{Symplectic \hcte and \phe disks}
\label{sympdiskpart}
First, we are going to show that the non-\spe \hcte of a compact
manifold with contact type boundary implies the existence of \phe
disk with a control on its area and its boundary:
\begin{thm}
\label{sympdisk}
Let $(M,\omega)$ be compact \spe manifold with contact type
boundary. If it is not \spe \hqe then, for any \pce \ste $J$
preserving the $\omega$-convexity (which means, preserving the
contact hyperplanes), there exists a $J$-\hle disk $f:\D \ra M$,
whose area is less than the capacity of $M$,
$\area(f(\D)) \leq \mu(M)$ and whose boundary is included in the
boundary of $M$ $f(\pp M) \subset \pp M$.
\end{thm}
\proof Let us choose $\mu \notin \mathcal{S}(\Sigma)$ such that
$FH^n(M,\mu) \ra H^{2n}(M,\pp M)$ is not surjective. Then there
exists $\eta >0$ such that $]\mu-\eta,\mu +\eta] \cap
\mathcal{S}(\Sigma) = \emptyset$.\\
Then according to \ref{rmqisom}, for any $\lambda > \mu$, $FH^n(M,\mu)
\ra H^{2n}(M,\pp M)$ is not surjective. Moreover, for $\lambda
>\mu$, it is possible to construct $\nu_\lambda$ and
$\epsilon_\lambda$ ($\epsilon_\lambda < \delta$) as small as
wanted (and thus $\displaystyle \lim_{\lambda \ra \infty}
\nu_\lambda =0$ $\displaystyle \lim_{\lambda \ra \infty}
\epsilon_\lambda =0$) such that $2 \nu_\lambda (
\mu-\eta)+\epsilon_\lambda < \eta$, $2 \nu_\lambda (
\mu+\eta)-\epsilon_\lambda < \eta$, and
$\epsilon_\lambda<\lambda \nu_\lambda$.\\
Then we consider the $\smooth$ Hamiltonian $H_\lambda$ on $\td{M}$
defined by:
\begin{equation}
\begin{cases}
(i) & -\epsilon_\lambda \textrm{ on } M \\
(ii) & h_\lambda(z) \textrm{ on } M_{[1,1+2 \nu_\lambda]},
\textrm{
  with } h_\lambda  \textrm{ a convex function on } [1,1+2 \nu_\lambda],\\
(iii) &  \lambda (z-(1+\nu_\lambda)) \textrm{ on } \Sigma
  \times [1+2\nu_\lambda,\infty).
\end{cases}
\end{equation}
Such a Hamiltonian exists because $\epsilon_\lambda +\lambda
\nu_\lambda < 2 \lambda \nu_\lambda$. It satisfies:
\begin{equation}
\begin{cases}
FH^*(K_\lambda,-\delta, \mu) \simeq FH^*(H_\lambda,-\delta, \mu) \\

FH^*(K_\lambda,-\delta, \delta) \simeq FH^*(H_\lambda,-\delta,
\delta).
\end{cases}
\end{equation}
Indeed there is a one-to-one correspondence between the
non-constant orbits of $H_\lambda$ and the closed characteristics
of $\Sigma$ whose action  $T \in \mathcal{S}(\Sigma)$ is equal to
$h_\lambda'(z)$ for some $z \in [1,1+2 \nu_\lambda]$ (which implies
$T<\lambda$). Then the associated orbit of $H_\lambda$ has action
$zT-h_\lambda(z)$. \\
And, by construction, if $T<\mu$ then $T<\mu-\eta$ and
$zT-h_\lambda(z) \leq (1+2\nu_\lambda)(\mu-\eta)+\epsilon_\lambda
<\mu$. Similarly, one checks that if $T>\mu$ then $zT-h(z) >\mu$.\\
In the same way, if $T>\delta$, then $T>T_0$, and, as $h_\lambda
(z)=-\epsilon_\lambda + \int_{1}^{z}h_\lambda '(z) \leq
-\epsilon_\lambda+(z-1)T$, we get $zT-h_\lambda (z) \geq T_0+
\epsilon_\lambda > \delta.$ Similarly one checks that if
$T<\delta$, then $zT-h_\lambda (z) < \delta.$\\

\par
Possibly restricting to one connected component of $M$,
$H^{2n}(M,\partial M) (\simeq H_0(M,
\partial M)) \simeq  \Z_2.[y_0]$, with $y_0$ a fixed point of $M$.
And let $\lambda
>\mu, \notin \mathcal{S}(\Sigma)$ be fixed. A $\smooth$ function
$f_\lambda$ is defined on $M$ by:
\begin{equation}
\label{deff}
\begin{cases}
(i)& f_\lambda \textrm{ is a Morse function and satisfies } -1
\leq
  f_\lambda \leq 0 \textrm{ on }  M_{1-\nu_\lambda}  \\

(ii)& f_\lambda \textrm{ has relative minimum only in }
y_0 \\
(iii)& f_\lambda(x,z)=z-(1-\nu_\lambda) \textrm{ on }
M_{[1-\nu_\lambda,1]}.
\end{cases}
\end{equation}

Then we consider, for small $\alpha >0$, the $\smooth$ Hamiltonian
$H_{\lambda, \alpha}$
\begin{equation}
H_{\lambda, \alpha} = \begin{cases}
(i) & -\epsilon_\lambda +\alpha \, f_{\lambda} \ {\rm on} \ M\\
(ii) &  g_{\lambda, \alpha}(z) \ {\rm on} \
M_{[1,1+2\nu_\lambda]},\ {\rm where} \
  g_{\lambda, \alpha}  \textrm{ is an arbitrary convex }\\
& \textrm{ function on } [1,1+2\nu_\lambda]  \textrm{ whose slope
   varies from } \alpha \ {\rm to}\  \lambda \\
&  \textrm{ and such that} \
  \lim_{\alpha =0} g_{\lambda, \alpha}= g_\lambda \ {\rm in} \ \smooth\\

(iii) & \lambda (z-1-\nu_\lambda) \textrm{ on } \Sigma \times
  [1+2\nu_\lambda,\infty).\end{cases}
\end{equation}

It is possible to perturb $J$ on $M_{[1-\nu_\lambda,1]}$ as little
as wanted into a generic structure $J_{\lambda, \alpha}$ so that
all trajectories linking an orbit inside of $M$ with an orbit in
$\Sigma \times
  [1,\infty) $ are admissible. This way we get $
  \lim_{\alpha \ra 0} J_{\lambda,\alpha} =J$ and $
  \lim_{\alpha \ra 0}  H_{\lambda, \alpha} =H_\lambda$ in
  $\smooth$. So, it follows that for $\alpha >0$ small enough

\begin{equation}
\begin{cases}
FH^*(H_{\lambda,\alpha},J_{\lambda,\alpha},-\delta, \mu) \simeq
  FH^*(H_\lambda,J, -\delta, \mu) \\

FH^*(H_{\lambda,\alpha} ,J_{\lambda, \alpha},-\delta, \delta)
\simeq FH^*(H_\lambda,J,-\delta, \delta).
\end{cases}
\end{equation}
Indeed, as $\delta$, $\mu \, \notin \mathcal{S}(\Sigma)$ , if
$\alpha$ is small enough, the Cerf diagram between $H_\lambda$ and
$H_{\lambda, \alpha}$ will be trivial (see \cite{viterbo}).

From the isomorphisms above, it is easy to check that, for
$\alpha$ small enough, the map
$FH^{2n}(H_{\lambda,\alpha},J_{\lambda,\alpha}, -\delta,\mu) \ra
FH^{2n}(H_{\lambda,\alpha},J_{\lambda,\alpha},-\delta, \delta)$ is
not onto. And so, since

 \begin {equation*}
\ra FH^{2n}(H_{\lambda,\alpha},J_{\lambda,\alpha},-\delta,\mu) \ra
FH^{2n}(H_{\lambda,\alpha},J_{\lambda,\alpha},-\delta,\delta) \ra
FH^{2n+1} (H_{\lambda,\alpha},J_{\lambda,\alpha},\delta,\mu) \ra
\end{equation*}
is an exact sequence, the map
\begin{equation}
\partial_{2n}~:FH^{2n}(H_{\lambda,\alpha},-J_{\lambda,\alpha},\delta,\delta)
\ra FH^{2n+1} (H_{\lambda,\alpha},J_{\lambda,\alpha},\delta,\mu)
\end{equation}
is not zero.\\

Moreover, by definition of $f_\lambda$ in (\ref{deff}),
$H^{2n}(H_{\lambda,\alpha},-\delta,\delta) \simeq
H^{2n}(M,\partial M) (\simeq H_0(M)) \simeq  \Z_2.[y_0]$ and so
$\partial_{2n}[y_0] \neq 0$. Then, there exists at least one orbit
$\gamma_{\lambda,\alpha} \in C^{2n+1}(H_{\lambda,\alpha},\delta,
\mu) $ such that $<\partial [y_0],\gamma_{\lambda,\alpha}> = 1$.
This implies $\delta \leq \mathcal{A}_H(\gamma_{\lambda, \alpha})
\leq \mu$ and so that $\gamma_{\lambda, \alpha}$ is one fo the
orbit associated to the closed characteristics of the boundary.\\

Therefore, there exists a Floer trajectory $u_{\lambda, \alpha}~:
\R \times \Ss^1 \lgra \td{M}$ between $\gamma_{\lambda, \alpha}$
and $y_0$ such that
\begin{equation*}
\begin{cases}(i) &\bar{\partial}_{J_{\lambda,\alpha}} u_{\lambda,
    \alpha}= \nabla H_{\lambda, \alpha} \\
(ii)&\lim_{s \ra + \infty} u_{\lambda, \alpha}(s,t) =y_0 \\
(iii)&\lim_{s \ra - \infty} u_{\lambda, \alpha}(s,.)
=\gamma_{\lambda, \alpha}\\
(iv)& E_H(u_{\lambda,\alpha}) = \mathcal{A}_H(\gamma_{\lambda,
\alpha}) - \mathcal{A}_H(y_0) \leq \mu
  + \epsilon_\lambda + \alpha. \end{cases}
\end{equation*}

Furthermore, by definition of $H_{\lambda, \alpha}$, necessarily
$\gamma_{\lambda, \alpha} \subset S_t$ with $t \in
[1,1+2\nu_\lambda]$. So, as the Floer trajectories can not be
tangent to any surfaces $\Sigma \times \{ z \}$  for $z\geq 1$, it
implies that $u_{\lambda,
  \alpha} (\R \times \Ss^1) \subset   M_{1+2\nu_\lambda}$.
\begin{rmk}
In order to get this trajectory, we first have to approximate
$H_{\lambda, \alpha}$  by non-autonomous Hamiltonians whose every
orbit is non-degenerate and then take the limit. For more details,
see \cite{her}.
\end{rmk}

Then, following \cite{her}, we can take the limit for $\alpha \ra
0$. As, $\displaystyle \lim_{\alpha \ra 0} H_{\lambda,
\alpha}=H_\lambda$ and $\displaystyle \lim_{\alpha \ra 0}
J_{\lambda,\alpha}=J$ in $\smooth$ we get trajectories $u_\lambda:
\R \times \Ss^1 \lgra M$ satisfying:
\begin{equation*}
\begin{cases}
(i) &\bar{\partial}_J u_\lambda= \nabla
H_\lambda \\
(ii) & \lim_{s \ra + \infty} u_\lambda =y_0\\
(iii) & \lim_{s \ra - \infty} u_\lambda =\gamma_\lambda,
{\textrm{ o{\`u} }}\  \gamma _\lambda\ \textrm{ est une
  orbite 1-periodique de }
  H_\lambda, \\ & \gamma_\lambda \subset
S_{t_\lambda} \ {\rm avec} \ t_\lambda \in [1, 1+2\nu_\lambda].\\
(iv) & E_{\rm{H}}(u_\lambda)  \leq \mu +
  \epsilon_\lambda.
\end{cases}
\end{equation*}

Let us notice that $\R \times \Ss^1$ can be identified with $\C^*$
by

\begin{equation}
\phi = \begin{cases} \R \times \Ss^1 & \lgra \C^* \\
                (s,t)& \lgra \exp(-(s+it))
\end{cases}
\end{equation}

Thus, $\td{u_\lambda}=u_\lambda \circ \phi$ is a map from $\C^*$
in $ M$, whose limit in $0$ is $y_0$ and which goes outside $M$
in $\infty$.\\
Then, there exists $z_\lambda \in (1-2\nu_\lambda, 1]$ (it can be
generically chosen as close as $1$ as needed) such that $(z \circ
\td{u_\lambda})^{-1} \{z_\lambda \}$ be a compact $1$-dimensional
submanifold $V_\lambda$ of $\C^*$. This $z_\lambda$ being fixed,
let $D_\lambda$ be the biggest topological disk of $\C$ containing
$0$ such that $\partial D_\lambda \subset V_\lambda$ et $0 \in
D_\lambda$. Such a disk exists because $\displaystyle \lim_{0}
u_\lambda =y_0$ and $z(y_0)\leq 1-\epsilon_0
<z_\lambda$.\\
As the Floer trajectories can not be tangent to $\Sigma \times \{
z \}$ for $z >1$, it follows that $\td{u_\lambda} ({D_\lambda}^*)
\subset M$.\\
Since $\triangledown H_\lambda= 0$ inside $M$ (and since
$\bar{\pp}_J u_\lambda = \dd \td{u_\lambda}(v)+ J \ \dd
\td{u_\lambda} (i\ v )$ with $v=-exp(-s-it)$),
$\bar{\partial}\td{u_\lambda}=0$ on ${D_\lambda}^*$ and  so
${u_\lambda}_{|{D_\lambda}^*}$ is $J$-\hl. Furthermore, this map
has a finite area. Indeed, being \ph, by denoting $Z_\lambda =
\phi^{-1}({D_\lambda}^*) \subset \R \times \Ss^1$, one writes:
\begin{equation*}
\displaystyle \area(\td{u_\lambda}({D_\lambda}^*))=
\int_{{D_\lambda}^*} \td{u_\lambda}^* \omega = \int_{Z_\lambda}
\left|\f{\pp (u_\lambda)}{\pp s}\right|^2 \leq E_{H}(u_\lambda)
\leq \mu +\epsilon_\lambda.
\end{equation*}
Moreover, as it has a finite limit on $0$, $u_\lambda$ can be
holomorphically extended to \phe function on $D_\lambda$.\\
Finally, one gets $\td{u_\lambda}$, \phe map from $\D_\lambda$,
with area less than $\mu + \epsilon_\lambda$ . Let $\theta$ be a
biholomorphism of $D_\lambda$ on the unit disk $\D$ that sends $0$
on $0$. Then $v_\lambda=u_\lambda \circ {\theta}^{-1}$ is a
$J$-holomorphic map from $\D$ into $M$ with $v_\lambda(0)=y_0$,
$v_\lambda(\partial \D) \subset M_{(1-2
  \nu_\lambda, 1]}$ and $\textrm{area}(v_\lambda)\leq \mu
+\epsilon_\lambda$ and thus for all $\lambda > \mu$. As
$\epsilon_\lambda$ and $\nu_\lambda$ go to zero when $\lambda$
goes to infinity, the theorem is proved.
\endproof  \\
As an immediate corollary:
\begin{cor}
Let $(M, \omega,J,\psi)$ be a subcritical Stein domain $\{ \psi
\leq a \}$. Then there exists a \phe disk $f:\D \ra M$ with area
less than $\mu(M)$ and with $f(\pp M) \subset \{\psi=a\}$.
\end{cor}

Now, let's consider an open $\omega$-convex manifolds $(M,
\omega)$. If it is non-\spe \hq, there exist an increasing
exhaustive sequence of $\omega$-convex domains $(M_j)$ for which
we have a precise control on the growth of the capacity (a
fortiori these capacities are finite and the $M_j$ are non-\spe
\hq). \\
Thus if $J$ is a compatible \pce \ste preserving the
$\omega$-convexity (more precisely preserving the contact
hyperplanes of each $M_j$), thanks to the theorem above, one gets
a sequence of \phe disks $f_j:\D \ra M_j$, $J$-holomorphic, on
which we have a control both on the area $\area(f_j) \leq
\mu(M_j)$ and on the position of the boundary $f_j(\pp \D) \subset
\pp M_j$. \\
In this context, thanks to some tools from complex analysis, we
could deduce from this the non-complex \hcte of $(M,J)$. This
reasoning will be specified and depending on the version of \spe
\hcte we consider, we will get different results. But first, let's
state the different needed results of (almost-)complex analysis.

\label{refforreasoning}

\section{(Almost-)complex tools}
\label{actoolspart}
We would like to get some results allowing to deduce from the
existence of a sequence of \phe disks with a control on their area
and their boundary, the non-complex \hct, and more precisely the
non-\aKe \hct.\\
A classical result is already available:

\begin{prop}
\label{lemmebgt} Let $(M,J)$ be a almost complex manifold and $g$
be a $J$-invariant Riemannian metric. If there exists
$J$-holomorphic maps $f_n:\D \lgra M$ such that
${area(f_n)}=O(d^2(f_n(0), f_n( \partial \D)))$, then
\begin{equation*}
\lim_{n \ra \infty} \sup_{z \in \D} |{f_n}'(z)|_{\rm{hyp}} =
\infty,
 \end{equation*}
\emph{i.e.} $M$ is not \aK-\hq.
\end{prop}
This proposition is an immediate consequence of the lemma:
\begin{lemma}
Let $(M,J)$ be a almost complex manifold and $g$ be a
$J$-invariant Riemannian metric. If $f: \D(1) \longrightarrow M$
is $J$-\hl, then $f$ is conformal and for every $r \in (0,1[$,
\begin{equation}
\max_{|z|\leq r} |f'(z)| \geq \frac{1}{r} ( \liminf_{|z| \nearrow
1} d_g(f(0),f(z))) - \frac{(\pi(1-r^2))^{\frac{1}{2}}}{2\pi r^2}
{(\mathrm{area}_g f)}^{\frac{1}{2}},
\end{equation}
\end{lemma}

In the particular case of $\C^n$, we can consider the standard
 Euclidean distance, that is to say the K{\"a}hler metric associated
 with $J_0$ and $\omega_0=dd^c \phi$ with $\phi=\frac{|z|^2}{4}$.
Thus, in this case, the hypothesis of proposition \ref{lemmebgt}
 reads as:  ${area(f_n)}=O(\min_{f_n(\partial \D)}\phi-\phi(f_n(0)))$.\\

It turns out that, for even more general Stein manifolds $(M,
\omega, \phi. J)$, the quantity $\min_{f_n(\partial
  \D)}\phi-\phi(f_n(0))$ is closely related to the invariant of Nevanlinna
$\tau(f_n)$: for Stein manifolds $\dsps
\tau(f)=\int_{0}^{2\pi}\phi \circ f(e^{i \theta}) d\theta \,-\,
2\pi \phi\circ f(0) $. In fact, $\tau(f)$ is defined for any
symplectic manifold by a more general expression but that appears
to coincide with this expression for the Stein manifolds.

That's why, we would like to get a similar result as the
proposition \ref{lemmebgt} comparing  ${\area(f_n)}$ with either
$\min_{f_n(\partial
  \D)}\phi-\phi(f_n(0))$ for Stein manifold, or $\tau(f_n)$ for general \aKe
manifolds. Thanks to the study of the invariant $\tau$, it is
possible to exhibit such a result, that appears to be stronger
than proposition \ref{lemmebgt} as explained beneath.

Let $(M,\omega, J)$ be a symplectic manifold provided with a
compatible almost complex structure  (that is an \aKe manifold).
Let us define:

\begin{defi}
If $f:\D(1) \longrightarrow M$ is a $J$-\hle map, we define
\begin{equation*}
\tau(f)=\int_{0}^{1} \frac{d\rho}{\rho} \int_{\D(\rho)} f^*
\omega.
\end{equation*}
\end{defi}

Then we show:
\begin{lemma}
\label{tau} $\tau(f)=\int_{\D} |f'(z)|^2
\ln\left(\frac{1}{|z|}\right)$ (with $|v|^2=\omega(v,Jv)$).
\end{lemma}
\proof Integrating by parts, one can check that:
\begin{equation*}
\tau(f)=\int_{0}^1 \ln\left(\frac{1}{\rho}\right)
\frac{\partial}{\partial
  \rho}\left(\int_{\D_{\rho}} f^* \omega\right) \rm{d} \rho.
\end{equation*}
As $f$ is $J$-holomorphic, it turns out that $f^* \omega=|f'(z)|^2
\rm{d}x \wedge \rm{d}y $. Thus,
\begin{equation*}
\tau(f)=\int_{0}^1 \ln\left(\frac{1}{\rho}\right)
\frac{\partial}{\partial
  \rho} \left(\int_0^\rho \left(\int_0^{2\pi} |f'(r e^{i\theta})|^2 r \rm{d}
  \theta\right)\rm{d}r\right),
\end{equation*}
what we wanted to show. \endproof  \\

As announced, in the particular case of Stein manifold, there is
another equivalent definition:

\begin{lemma}
\label{taupsi} If $\omega=dd^c \psi$, then
\begin{equation}
\label{eq1} \tau(f)=\int_{0}^{2\pi}\psi \circ f(e^{i \theta})
d\theta \,-\, 2\pi \psi\circ f(0).
\end{equation}
\end{lemma}
\proof Let us note $g=\psi \circ f$. Then, as $f$ is $J$-\hl,
$\tau(f)=\int_0^1 \left(\int_{\D(\rho)}dd^c g\right)
\frac{\rm{d}\rho}{\rho}$. Applying Stokes theorem, we get
\begin{equation*}
\tau(f)=\int_0^1 \left(\int_{\Ss(\rho)} d^c
g\right)\frac{\rm{d}\rho}{\rho} = \int_0^1 \left(\int_0^{2\pi}
\rm{d}g_{|\rho e^{i\theta}}(e^{i\theta}) \rm{d}\theta\right)
\rm{d}\rho =\int_0^1  \frac{\partial}{\partial
  \rho} \left(\int_0^{2\pi} g(\rho  e^{i\theta})\rm{d}\theta \right)\rm{d}\rho.
\end{equation*}
\endproof \\
\indent
 The lemma \ref{tau} provides us with the wished result
 linking the existence of \phe disks, with control on both the area
 and the boundary position, and the \aKe \hct:

\begin{prop}
\label{proptau} If there exist $J$-\hle maps $f_n~: \D \lgra M$
such that
 $ \area(f_n) = O(\tau(f_n))$, then:
\begin{equation*}
\lim_{n \ra \infty} \sup_{z \in \D} |{f_n}'(z)|_{\rm{hyp}} =
\infty.
 \end{equation*}

\par Thus $(M,\omega,J)$ is not \aKe \hqe (thus, in the several compact cases
mentioned in section \ref{complexhyppart}, $M$ is not complex \hq).
\end{prop}

\begin{cor}
Let $(M,\omega, J, \psi)$ s Stein manifold. If there exist
$J$-\hle maps $f_n:\D \ra M$ with $f_n(\pp \D) \subset
\{\psi=a_n\}$, $f_n(0)=x_0$ and $\area(f_n)=O(a_n)$, then $M$ is
not \aKe-\hq.
\end{cor}
This is a immediate corollary from the proposition above since,
according to the lemma \ref{taupsi}, under these hypothesis
$\tau(f_n)=2 \pi(a_n-\psi(x_0))$.

\proof (of proposition \ref{proptau}) In sight of lemma \ref{tau}
we have,
\begin{equation*}
\tau(f_n)= \int_0^{2\pi} \int_0^1 |{f_n}'(r {\rm{e}}^{i
\theta})|^2 \ln \left( \f{1}{r} \right) r {\rm{d}}r
{\rm{d}}\theta.
\end{equation*}

Thus, by cutting the inside integral between $0$ and $1$ in $R$
with $0<R<1$~:
\begin{equation}
\int_0^1 |{f_n}'(r {\rm{e}}^{i \theta})|^2 \ln \left( \f{1}{r}
\right) r dr \leq A \,R \, \sup_{r \leq R} |{f_n}'(r \rm{e}^{i
\theta})|^2 + ln \left( \f{1}{R} \right) \int_R^1 |{f_n}'(r
\rm{e}^{i \theta})|^2 r {\rm{d}}r \rm{d}\theta,
\end{equation}
with $\displaystyle A=\max_{[0,1]} \ln \left( \f{1}{r} \right) r$.
Integrating with respect to $\theta$, we get $\tau(f_n) \leq A \,R
\, \sup_{\D_R} |f_n'(z)|^2 \, + area(f) \ln \left( \f{1}{R}
\right)$. Thus,
\begin{equation}
 \sup_{\D_R} |f_n'(z)|^2 \geq \f{\tau(f_n)}{A}\left( 1 - \ln \left( \f{1}{R}
\right) \f{area(f_n)}{\tau(f_n)} \right).
\end{equation}
By hypothesis, there exists $B>0$ such that
$\f{area(f_n)}{\tau(f_n)}
  \leq B$.
Therefore, there exists $R<1$ such that $\displaystyle  \forall n,
\,\sup_{\D_R} |f_n'(z)|^2 \geq \f{\tau(f_n)}{2A} \, \ra \infty$.
Thus we get the wanted result (because $\sup_{\D_R}
|f_n'(z)|_{hyp}^2 \geq (1-R^2) \sup_{\D_R} |f_n'(z)|^2$).
\endproof \\
Let us point out that if $g=\omega(.,J.)$ is the \aKe metric
associated to $\omega$ and  $J$, then for every $J$-\hle map $f$,
$ \tau(f) \geq4 \pi \, d_g^2(f(0), f(\pp \D))$. Indeed for every
$\theta$, \ba* d_g^2(f(0), f(\pp \D)) & \leq & \left( \int_0^1
|f'(r {\rm e}^{i\theta}|
  dr\right)^2 \\ & \leq&  \left(\int_0^1 |{f_n}'(r \rm{e}^{i \theta})|^2 \ln
\left( \f{1}{r} \right) r dr \right) \, \left(\int_0^1 \f{1}{ \ln
\left( \f{1}{r} \right) r} dr \right), \ea*
with the last integral
being equal to $\f{1}{2}$. Thus, integrating the inequality with
respect to $\theta$, we get the announced inequality.
That's why proposition
\ref{proptau} appears to be stronger than proposition
\ref{lemmebgt}: the hypothesis is weaker and can be used even for
manifolds with bounded diameter (see examples of polarizations for
example).\\

\par
Let us mention another proposition that allows, strengthening the
hypothesis, to show directly the non-Kobayashi hyperbolicity (the
previous proposition implies the non-Kobayashi \hcte only in the
compact case):
\begin{prop}
\label{proptau2} Let  $(M,\omega, J, \psi)$ be a Stein manifold.
If there exist $J$-\hle maps $f_n~: \D \lgra M$ such that
$\displaystyle \lim_{n \ra
  \infty} \frac{\tau(f_n)}{(\area(f_n))} = \infty$ (that is
$\area(f_n)=o(\tau(f_n))$) and $f_n(0)=x_0$, then:

\begin{equation*}
\forall \, r>0, \, \, \lim_{n \ra \infty} \sup_{z \in \D_r}
|{f_n}'(z)|_{\rm{hyp}} = \infty,
 \end{equation*}
and $M$ is not Kobayashi-\hq.
\end{prop}

\proof Let us note $C_r=[0,r] \times [0, 2\pi]$, $A_r=[r,1] \times
[0, 2 \pi]$ and
$\tau_r(f_n)=\int_{C_r} |f_n'|^2 \ln\left(\frac{1}{t}\right) t dt d\theta.$ \\
Then, $\tau(f_n)=\tau_r(f_n) + \int_{A_r} |f_n'|^2
\ln\left(\frac{1}{t}\right) t dt d\theta.$. Thus,
\begin{equation}
\label{minoration} \tau_r (f_n) \geq \tau(f_n) -\ln
\left(\frac{1}{r}\right) \area(f_n(A_r)) \geq \tau(f_n)
\left(1-\ln\left(\frac{1}{r}\right)
\frac{\area(f_n)}{(\tau(f_n))}\right).
\end{equation}

As in the proof of the previous proposition, one can check that
there exists a constant $A$ such that for all $n$, $\tau_r(f_n)
\leq A \sup_{\D_r} |f'_n|_{\rm{hyp}}$.\\
 Therefore, in sight of
(\ref{minoration}), we get the first announced
result.\\
In order to deduce that $M$ is not Kobayashi-\hq, we need the
lemma:

\begin{lemma}
Let $V_0$ be a relatively compact open set containing $x_0$. Then,
$\forall \, r >0, \, \exists n_0 / \, \forall n >n_0, \, f_n(\D_r)
\not\subset V_0.$
\end{lemma}
Indeed, similarly to lemma \ref{taupsi}, one can check that
$\int_0 ^{2\pi} (\psi \circ f_n)(r e^{i \theta}) d\theta -2\pi
\psi \circ f_n(0) =
\left(-\ln\left(\frac{1}{r}\right)\right)\left(\int_{D_r} f_n^*
\omega\right) + \tau_r(f_n).$ Since $\psi$ is bounded on $V_0$,
there  exists a constant $A$ (depending only on $V_0$) such that,
if  $f_n(\D_r) \subset V_0$, then
\begin{equation*}
\ln\left(\frac{1}{r}\right) \geq \frac{1}{\area(f_n(\D(r)))}
\left(\tau_r(f_n) -A\right).
\end{equation*}
In sight of the inequality (\ref{minoration}), we get
$\ln\left(\frac{1}{r}\right) \geq \frac{1}{\area(f_n)}
\left(\tau(f_n) -A\right)$, which goes to infinity. Thus the lemma
is proved. \endproof \\
 Then, $\forall m$, $\exists n_m$ $ /$  $
f_{n_m}(D_\frac{1}{m}) \not\subset V_0$. Therefore there exists
$z_m \in D_{\frac{1}{m}}$ such that $p_m=f_{n_m}(z_m) \in \partial
V_0$. If we note $d_K$ the Kobayashi pseudo-distance, these points
satisfy $d_K(x_0, p_m) \leq d_{\rm{hyp}}(0,z_m)$, whose limit is
$0$ when $m \ra \infty$. We can extract a subsequence of  $(p_m)$
admitting a limit $p \in
\partial V_0$ (thus $p \neq x_0$). Then $d_K(p, x_0)=0$ and $M$ is not
Kobayashi-\hq. \endproof \\

\section{First links between \spe and complex hyperbolicities}
\subsection{Main theorems}
\label{maintheoremssection}
The reasoning explained at the end of the section
\ref{refforreasoning} and the complex tools of the previous
section provide us with some wished theorems linking \spe and
complex hyperbolicities. First in the most general case,

\begin{thm}
\label{resdist}
Let $(M^{2n},\omega)$ be an open $\omega$-convex \spe manifold
satisifying the non-\sp-\hcte assumption:
there exists an exhaustive increasing
sequence $(M_j)$ with $\mu(M_j)= O(d_g^2(\pp M_j, x_0))$. Then for any
(uniformly for the fixed metric $g$) compatible almost-complex \ste
 preserving the
$\omega$-convexity (of the $M_j$), $(M,J)$ is not \aK-\hq.
\end{thm}
The contrapositivity reads:
\begin{thm}
Let $(M,\omega)$ be a $\omega$-convex \spe manifold. If there exists $J$
a compatible \pce \ste such that $(M,J)$ is \aK-\hq, then (if $g$ denotes
the \aKe metric associated with $\omega$ and $J$), $(M,\omega)$ is
symplectic \hqe: for any increasing exhaustive sequence of $J$-convex
domains $(M_j)$,
\beq*
\f{\mu(M_j)}{d_g^2(\pp M_j, x_0)} \ra \infty.
\eeq*
\end{thm}
\proof (of theorem \ref{resdist}) Let $(M_j)$ be  a sequence of
$\omega$-convex domains such that $\mu(M_j)= O(d_g^2(\pp M_j,
x_0))$. As explained, if $J$ is a compatible \pce \ste preserving
the $\omega$-convexity of each $M_j$, thanks to theorem
\ref{sympdisk}, for each $j$, there exists a $J$-\hle disk $f_j:\D
\ra M$ with $f(0)=x_0$. $f_j(\pp \D) \subset  \pp M_j$ and
$\area(f_j) \leq \mu(M_j) = O(d_g^2(\pp M_j, x_0))=O(d_g^2(f_j(0),
f_j(\pp \D))$. Then the proposition
\ref{lemmebgt} provides us with the theorem. \endproof\\
For closed manifolds we can immediately deduce (in this case all
the complex hyperbolicities are equivalent):
\begin{thm}
Let $(M, \omega)$ be a closed \spe $\omega$-convex manifold. If it is non
\sp-\hq, \emph{i.e.} if there exists an exhaustive increasing
sequence $(M_j)$ with $\dsps \td{M}=\cup_i M_j$ such that
$\mu(M_j)= O(d_g(x_0,\pp M_j)^2)$, then for any compatible \pce
\ste that respects the $\omega$-convexity (which means that
respects the contact hyperplanes of each $\pp M_j$), $(M,J)$ is
not complex \hq.
\end{thm}

In the particular case of Stein manifolds, the considered
$\omega$-convex domains are the $M_a=\{ \psi \leq a \}$. So by
applying the same reasoning than above and using the proposition
\ref{proptau}, we prove:
\begin{thm}
\label{res2} Let $(M, J, \omega, \psi)$ be a Stein manifold and
$M_a=\{ \psi \leq a\}$. If $M$ is not (S-)\spe \hq, more precisely
if there exists an exhaustive sequence $(a_n)$ such that
$\mu(M_{a_n}) = O(a_n)$, then $M$ is not \aK-\hq.
\par
Moreover, if $M$ satisfies the stronger assumption of non-\sp-\hct:
$\mu(M_{a_n})= o(a_n)$, then $M$ is not
Kobayashi-\hq.
\end{thm}
Let's notice than the non-\aK-hyperbolicity implies the
non-complex hyperbolicity only under compactness assumptions. So,
the second part of this theorem happens to be very useful in the
non compact case. The proof is similar as the previous one:\\
\proof According to the hypothesis and theorem \ref{sympdisk},
there exists a sequence of bigger and bigger \phe disks: $f_a:\D
\ra M$ whose origin $f_a(0)=x_0$ is fixed, whose boundary position
is controlled by $f_a(\pp \D) \subset \pp M_a =\{ \psi=a \}$, and
whose area is controlled by $\area(f_a) \leq \mu(M_a)$. Thus
$\tau(f_a) = 2 \pi \, a +c $ (with $c=-2 \i \psi(x_0)$
constant)and then the proposition \ref{proptau} and \ref{proptau2}
provide us with the wished result.
\endproof  \\
This theorem also reads:
\begin{thm}
Let  $(M, J, \omega, \psi)$ a Stein manifold and  $M_a=\{ \psi
\leq a\}$.  If $M$ is \aK-\hq, then $(M,\omega)$ is  \spe \hqe:
\beq* \f{\mu(M_a)}{a} \ra \infty \eeq* \\
If $M$ is Kobayashi-\hqe then there exists a constant $C>0$ such
that \beq* \mu(M_{a}) \geq C a \eeq*
\end{thm}
Furthermore, according to the particular estimate we got for the
growth of the capacity in the special case of Stein manifold,
theorem \ref{res2} has as a corollary:
\begin{cor}
\label{reswithalpha} Let $(M^{2n}, J, \omega)$ be a Stein
manifold, complete at $-\infty$, (\emph{i.e.} such that the
gradient flow of $\psi$ is complete at $-\infty$). If $M$ is not
strongly \sp-\hqe (for example if $M$ is subcritical) and if there
exist a constant $\alpha>0$ and a compact $K \subset M$ such that
$\frac{||d\psi||^2}{\psi} \geq \alpha$ on $W \exc K$, then:
\begin{description}
\item[\quad $(i)$] if $\alpha \leq 1$,
  then $M$ is not  \aK-\hq.
\item[\quad $(ii)$] if $\alpha < 1$, then $M$ is not
Kobayashi-\hq.
\end{description}
\end{cor}
These result allow us to get some results about the \hcte of Stein
manifolds. For example it applies to the case of manifolds $W=M
\exc H$ with $M$ a closed complex manifold and $H$ some
hypersurface. Let's explain this.\\

\subsection{Application to the manifolds $W \exc H$}
\label{applicationw}
Let's consider a closed \Kae manifold $(M,\omega_0,J_0)$. If
$\f{1}{2\pi} [\omega_0]$ is an entire cohomology class
(\emph{resp.} rational), there exists a pre-quantization of
$(M,\omega_0)$ (\emph{resp.} $(M, k \omega_0)$), which means a
\hle line bundle over $M$ provided with a hermitian metric and a
connection $(\mcl{L}, |.|, \nabla)$ with curvature form $\omega_0$
(\emph{resp.} $k \omega_0$, in which case, $k \omega_0$ will be
denoted itself by $\omega_0$ in the following). For any nonzero
\hle section $s:M \ra \mcl{L}$, one can define the complex
hypersurface $H=\{s=0\}$ and the manifold $W=M\exc H$. Then, the
function $\phi=(-log |s|^2)$ is pluri-sub-harmonic on $W$ and
satisfies $\dd \dd^c \phi= \omega_0$. Thus $(W, \omega_0, J_0,
\phi)$ is a Stein manifold (non complete).\\
So the above results can be applied and this particular case, the
involved quantity $\frac{||d\phi||^2}{\phi}$ can be more precisely
estimated. Indeed, $|d\phi|=\f{1}{|s|^2} |\dd|s|^2|$. Thus, \be
\label{inegverifnormesigma} \frac{||d\phi||^2}{\phi}=\f{|\pp
|s|^2|^2}{|s|^4 \ |\log|s|^2|}. \ee

Let's generalize the notion of transverse section defining a
section $s$ to be {\bf pseudo-transverse to the null section} if
there exist a constant $C>0$, a constant $0 \geq \alpha <1$,
and a neighborhood $V$ of $\Sigma$
such that $|d|s|^2|> C|s|^{2 \alpha}$ on $V$. Thus according to the
inequality (\ref{inegverifnormesigma}) and proposition \ref{resmu}, we have
\begin{lemma}
\label{cappseudotrans1}
Let $(\mcl{L},|.|, \nabla)$ a hermitian \hle line
bundle over the \Kae manifold $(M, \omega_0, J_0)$ with curvature
form $\omega_0$. For $s$ any nonzero \hle section of $\mcl{L}$
which is pseudo-transverse to the null-section, if $W=M\exc \Sigma$,
with $\Sigma=\{s=0\}$, is not strongly \spe \hq, then for any
$\alpha$, there exists a constant $c>0$ such that $\mu(\{\phi \leq a\},
\omega_0) \leq c \ a^\alpha$.
\end{lemma}
And the corollary \ref{reswithalpha} reads:
\begin{prop}
\label{propex} Let $(\mcl{L},|.|, \nabla)$ a hermitian \hle line
bundle over the \Kae manifold $(M, \omega_0, J_0)$ with curvature
form $\omega_0$. For $s$ any nonzero \hle section of $\mcl{L}$
which is pseudo-transverse to the null-section (in particular for
any transverse section) let's consider $\Sigma=\{s=0\}$ et
$W=M\exc \Sigma$. If $W$ is not strongly \spe \hq, then
\begin{enumerate}
\item $(W, \omega_0, J_0)$ is not \aK-\hq. So $(W, J_0)$ is not
hyperbolically embedded in $M$. And either $W$ or $\Sigma$ is not
  Brody-\hqe (and so $M$ is not\hq).
\item  $(W,J_0)$ is not Kobayashi-\hq.
\end{enumerate}
\end{prop}
In particular this can be applied for the polarizations with the
definition introduced by Biran and Cielibak \cite{biran}:
$\mathcal{P}=(M,\omega, J, \Sigma)$ is a polarization of the \Kae
manifold $(M,\omega,J)$ if $[\omega] \in H^2(M, \Z)$ and if
$\Sigma$ is reduced complex hypersurface, Poincare dual of $k
[\omega]$, where
$k$ is called the degree of the polarization.\\
Then the reasoning above can be applied to the bundle
$\mathcal{L}=\mathcal{O}_M(\Sigma)$ after choosing a \hle section
$s$ of $\mathcal{L}$ defining $\Sigma$.\\
A polarization is called sub-critical if the Stein manifold $W=M
\exc \Sigma$ is a sub-critical Stein manifold. Then $W$ is not
strongly \spe \hq. And, as a corollary of the proposition above,
one gets:
\begin{cor}
If the closed \Kae manifold $(M,\omega,J)$ can be provided with a
sub-critical polarization $M,\omega, J, \Sigma)$ defined by a
pseudo-transverse section $s$ of the bundle
$\mathcal{O}_M(\Sigma)$, then $(M,J)$ is not Brody-\hq. More
precisely, either $W=M \exc \Sigma$ or $\Sigma$ is not Brody-\hq.
Besides, $(W,J)$ is not Kobayashi-\hq.
\end{cor}
The proposition \ref{propex} allows to exhibit some more examples
of non-complex \hqe manifolds.\\
It notably applies to the case of the manifold $\PP^n \C$, provided
with its standard \st s ($\omega$ and $J$), minus $k$ hyperplanes
in general position, with $k \leq n$. Indeed such a hypersurface
is pseudo-transverse: the coordinates system can be chose so that
the $k$ hyperplanes are defined by $z_j=0$ for $0 \leq j \leq
k-1$. then the associated pluri-sub-harmonic function is $\phi=-
\log (|\sigma|^2)$ with

\beq* \dsps |\sigma| = \f{\Pi_{j=0}^{k-1} |z_j|}{|z|^{k}}=0.
\eeq* And $\dd \dd^c \phi =k \omega$. \\
For any $j=0 \dots k-1$, by applying $\dd \phi$ to the vector $(0,
\dots, z_j, \dots, 0)$, one notices that on the neighborhood
$V_j=\{ |z_j|^2 \leq \f{1}{2k} \}$ of the hyperplane $\{z_j=0\}$,
$|d\phi|^2 \geq \f{c}{ |z_j|^2}$. So, on the neighborhood of
$\Sigma$, $V=\displaystyle \cup_{j=0}^{k-1} V_j$ (\emph{i.e.}
outside a compact set of $W$),
\begin{equation*}
|d\phi|^2 \geq c \max_{j=0 \ldots k-1} \f{1}{|z_j|^2} \geq c
\f{1}{|\sigma|^{\f{1}{k}}}.
\end{equation*}
Finally on $V$,  \beq* \frac{||d\phi||^2}{\phi} \geq c
\f{1}{|\sigma|^{\f{1}{k}}\ |\log(|\sigma|^2)} \ra \infty \eeq*
when $|\sigma| \ra 0$ and so the section $\sigma$ is
pseudo-transverse. Thus the proposition \ref{propex} allows to get
the (already known) result: $\PP^n \C \exc \{k \textrm{ hyperplanes} \}$ is
not complex-\hq.

\section{Deepening: stability of non-complex-\hct}
\label{deepeningpart}
In the previous part, we have shown that the non-symplectic \hcte
implies the non-complex \hcte of any \pce \ste "that preserves the
$\omega$-convexity". However because of this hypothesis, the
result is not completely satisfactory and we would like to get rid
of this assumption. One of the main motivation is to generalize
Bangert's theorem \ref{bangertthm} for some more general manifolds
by proving the
non-complex \hcte for any \pce \ste compatible with the
non-\sp-\hqe \spe \st. Or at least, we would like to prove the
stability of the non-complex \hcte by small deformations of the
\pce \ste inside the sets of compatible \st s.\\
The previous study already provides us with a result
of stability for non-complex-\hct, straightforward consequence of
theorem \ref{resdist}.

\begin{prop}
\label{primsstab}
Let $(M, \omega)$ be an open $\omega$-convex \spe manifold and
$J_0$ be a compatible \pce \ste such that $M$ can be written as an increasing
union of relatively compact $J_0$-convex domains $M_j$. Let's denote
 $|.|$ and $d$ the norm and distance associated to $g_0=\omega(.,J_0.)$.\\
If $(M, \omega)$ is not $g_0$-\spe \hq,then $(M, \omega, J)$ is
not \aK-\hqe for any $J$ in an open $\mcl{C}^1$-neighborhood of $J_0$,
more precisely for any compatible $J$ such that $g=\omega(.,J.)$ and
$g_0$ are equivalent metric ($\mcl{C}^0$-condition) and that the
$M_j$ are $J$-convex (which is the case if there exists some
$\alpha<1$ such that $\dsps |\dd \dd^c_{J-J_0} \phi| < \alpha
\inf_{v \in \zeta_x, |v|=1} \dd \dd^c_{J_0} \phi (v,Jv)$, where
$\phi$ defines $\pp M_J=\{ \phi=a_j \}$ ).
\end{prop}

\proof Let the hypersurface $\pp M_j$ be defined by $\{\phi
=a_j \}$  with $\phi$ fonction on $M$ and $a_j$ some constants.
Since they are $J_0$-convex, $\dd \dd^c_{J_0} \phi $ is strictly positive
on the complex hyperplanes of $T\pp M_j$. If there exists some
$\alpha <1$ such that
$\dsps |\dd \dd^c_{J-J_0} \phi| < \alpha \inf_{v \in \zeta_x,
|v|=1} \dd \dd^c_{J_0} \phi (v,Jv)$, then
$\dd \dd^c_J(v,Jv)=\dd \dd^c_{J_0} \phi
(v,Jv) +\dd \dd^c_{J-J_0} \phi (v,Jv) \geq \dd \dd^c_{J_0}
\phi(v,Jv) -\alpha \inf_{v \in \zeta_x, |v|=1} \dd \dd^c_{J_0}
\phi (v,Jv) |v|^2$. Moreover if $J$ is close enough to $J_0$,
the complex hyperplanes for $J$ are close to the ones for $J_0$.
So finally there $\alpha < \alpha" <1$ such that on these hyperplanes
$\dd \dd^c_{J_0} \phi (v,Jv)\geq \alpha"\ \inf_{v \in \zeta_x,
|v|=1} \dd \dd^c_{J_0} \phi (v,Jv) |v|^2$. Thus the $\pp M_j$ are
$J$-convex. Besides, if there exists
$\alpha' <1$ such that $|J-J_0|< \alpha'$ then the metrics
$g$ et $g_0$ are equivalent and then $\mu(M_j)
=O(d(x_0,M_j)^2= O(d_g(x_0,M_j)^2)$. Thus theorem \ref{resdist}
can be applied to $(M,\omega,J)$ which concludes the proof.
\endproof\\
However the condition required on $J$ is quite strong (thus the open set
of non-complex-\hqe \st s is quite small). And we'd like to improve this
result.

\subsection{Main results}
\label{enouncethebig}
The idea will be to fundamentally use the contractibily of the set
of compatible \pce \st s. A compatible \pce \ste being fixed,
thanks to this property, it could be deformed in \pce \st s
preserving the $\omega$-convexity at infinity. Thus we get a
family of \phe curves for these \pce \st s. Then we will need some
isoperimetric assumptions to cut these curves and get a \phe curve
for the initial fixed \pce \st. These isoperimetric assumptions are
in particular satisfied by the
manifold whose geometry is well-bounded in the meaning of the next
definitions. \\
\par
But first, let's recall that a riemannian manifold is said to be with
conical ends, if he can be decomposed as
$\dsps M=M_0 \cup \cup_j W_j$, with $M_0$ compact with boundary
$\pp M_0 =\cup_j S_j$, and $W_j$ some conical ends diffeomorphic to
$S_j \times [1, \infty[$, $M$ being built by gluing the $W_j$ on $M_0$
along $S_j$. Moreover, there exists on each $W_j$ a (Busemann or distance)
function $\rho_j$ such that $|\f{\pp}{\pp \rho_j}|=1$ and $S_j \times \{t\}
=\rho_j^{-1}(\{t\})$. Then the radial curvature is defined as the sectional
curvature restricted to the plans containing $\f{\pp}{\pp \rho}$. The
particular case of the $W_j$ being hermitian is studied in \cite{greenewu}.
It is notably proved that if the radial curvature is bounded from above
on $\{\rho =s \}$ by a function $K(s)$ satisfying $\int K<1-\mu <1$ then
$Hess \rho \geq \f{\mu}{\rho} H$ with $H=g_0 - d\rho \otimes d\rho$ (and
so $\dd \dd^c \rho^2 \geq 4 \mu g$, \emph{i.e.} for any $v$,
$\dd \dd^c \rho^2(v,Jv) \geq 4 \mu |v|^2$). This is this assumption on the
Hessian that will be needed in our study. However, according to this
result, it can be intuitively considered as, in some way, an assumption
of bounded radial curvature.\\
As a remark, let's notice that the assumption on the radial curvature and
especially the boundary of integration of $\int K$ need to be specified.
Writing this, one implies that the ends $W_j$ can be extended in a pole
manifold $\bar{W_j}$ on which $\int_0^\infty K(s)<1$ (see \cite{greenewu}
for more details).\\
Let's now introduce all the other assumptions of bounded geometry that will
occur in our study:

\begin{defi}
A Stein manifold $(M,\omega,J, \psi)$ is said to have
\begin{enumerate}
\item a well-bounded Hessian if there exists a constant $A>0$ such
that $\forall v \in \{ d\psi=0\}$ $Hess\,\psi(v,v) \geq A
\f{|\nabla\psi|^2}{\psi} |v|^2$

\item a $\delta$-bounded derivative (with $\delta \geq 0$) if the
canonical function $\beta$ (defined as $\dsps \sup_{x,
\sqrt{\psi}=s}\f{|\nabla\psi|^2}{\psi}$) is bounded from above:
$\beta(s)= O(s^\delta)$.

 \item a radial curvature well-bounded if it can be written as a manifold
 with conical ends such that \begin{itemize}
\item either $Hess \rho \geq \f{\mu}{\rho} H$ with
$H=g_0 - d\rho \otimes d\rho$ (\emph{i.e.} $Hess \rho \geq g_0$ on $TS$)
\item or there exists a
function $K\geq 0$ such that the radial curvature on $\{ \rho=s
\}$ is bounded from above by: $curv \leq K(s)$, and $\int s K(s)
<1$ (which implies the previous assumption)
\end{itemize}
\item  a $I_\delta$-bounded geometry if its Hessian is
well-bounded and its derivative is $\delta$-bounded

\item a $II_\delta$-bounded geometry if its radial curvature is
well-bounded and its derivative is $\delta$-bounded

\end{enumerate}
\end{defi}

\begin{exam} In the case of $(\C^n, \omega_0, J_0, g_0, \psi)$,
$\psi=\f{|z|^2}{4}$, $\f{|\nabla\psi|^2}{\psi}=1$ so its
derivative is well-bounded. $Hess \psi (v,v) =\f{1}{2} |v|^2$, so
its Hessian is (even strongly) well-bounded. Moreover its radial
curvature is null so a fortiori well-bounded. And at the end
$\C^n$ geometry is both (even strongly) $I_0$-bounded and
$II_0$-bounded.
\end{exam}
In fact these assumptions of Hessian (respectively radial curvature) bounded
curvature happen to be useful in our study because they imply that the
projection on the sets $M_a$ (respectively $B(a)$) are contracting enough,
which is what we need to cut the \ph s curves. That's why, in all the
following theorems and corollaries, these hypothesis could in fact be
replaced by one on the contracting properties of the projection (as it will
appear in the proof of lemma \ref{implicationlemma}: see proof of
corollary \ref{corminriem} and proposition \ref{minimstein}).\\
\par

Let's now state the results:

\begin{thm}
\label{thebigone1} Let $(M, \omega, \psi, J_0)$ be a Stein
manifold. If it satisfies one of the following properties:

\begin{enumerate}
\item \label{one} \begin{itemize}

\item For $Q>1$ there exists a constant $C=C_Q>0$ such that if a
$Q$-minimizing surface $S$ satisfies $\pp S \subset B(x_0,C
\sqrt{\mu(M_{a+1})})$ and $\area (S) \leq \mu(M_{a+1})$ then $S
\subset M_a$, and

\item $(M,\omega)$ is not strongly \sp-\hq.
\end{itemize}

\item \label{three} Or, \begin{itemize}

\item For any $Q>1$ there exists a constant $C>0$ such that if
a $Q$-minimizing surface $S$ satisfies $\pp S \subset M_{C \, a}$
and $\area(S) \leq \mu(M_{a+1})$ then $S \subset M_a$, and

\item  $M$ satisfies the non-\sp-\hcte assumption: there exists an
exhaustive sequence $(a_n)$ such that
\begin{equation}
\label{n2} \mu(M_{a_n+1})=O\left(\f{a_n}{\beta(a_n)}\right)
\end{equation}
\end{itemize}
\end{enumerate}
Then for any $J$ uniformly compatible with $\omega$, $(M,
\omega,J)$ is not \aK-\hq. Thus if $(M, \omega)$ possesses a
compact quotient $(W, \omega_0)$, then for any $J$ compatible with
$\omega_0$, $(M, J)$ is not Brody-\hq.
\end{thm}
One can deduce from this theorem another version by enouncing the
assumptions in terms of bounded geometry rather than in terms of
isoperimetric properties (that they imply).

\begin{thm}
\label{thebigone2} Let $(M, \omega, \psi, J_0)$ be a Stein
manifold. If it satisfies one of the following properties:

\begin{enumerate}
\item \label{six} The geometry of $M$ is $II_\delta$-bounded. And
$M$ satisfies the non-\sp-\hcte assumption:
\begin{equation}
\label{autreeq} \mu(M_{a})=O(a^{(1-\delta)\f{m}{2}}),
\end{equation}
with $m$ an isoperimetric value of $M$.

\item \label{five} The geometry of $M$ is $I_\delta$-bounded. And
$M$ satisfies the non-\sp-\hcte assumption:

\beq* \mu(M_{a})=O\left( a^{\min(1-\delta,\, (1+\delta)
\f{m}{2})}\right) \eeq*
\end{enumerate}
Then, for any $J$ uniformly compatible with $\omega$, $(M,
\omega,J)$ is not \aK-\hq. Thus if $(M, \omega)$ possesses a
compact quotient $(W, \omega_0)$, then for any $J$ compatible with
$\omega_0$, $(M, J)$ is not Brody-\hq.
\end{thm}

If the manifold satisfies some less strong non-\sp-\hcte assumptions,
then we
couldn't prove anymore the non-complex-\hcte for any compatible \pce \ste
anymore. However, we can prove the non-complex-\hcte for any
compatible \pce \ste in an open neighborhood of the standard one.
\begin{thm}
\label{thebigone3} Let $(M, \omega, \psi, J_0)$ be a Stein
manifold. If it satisfies one of the following properties:

\begin{enumerate}
\item \label{two} For any $Q>1$, there exists a constant $C>0$
such that if a $Q$-minimizing surface $S$ satisfies  $\pp S
\subset M_{C \, a}$ and $\area(S) \leq \mu(M_{a+1})$ then $S
\subset M_a$. Moreover, $M$ is not (S-)symplectic-\hq: there
exists an exhaustive sequence $(a_n)$ such that $\mu(M_{a_n})
=O(a_n)$.

\item \label{four} The geometry of $M$ is $I_\delta$-bounded.
Moreover $M$ satisfies the non-\sp-\hcte assumption:
\beq* \mu(M_{a})=O\left( a^{\min(1,\, (1+\delta) \f{m}{2})}\right)
\eeq*
\end{enumerate}

Then, $(M, \omega, J_0)$ is not \aKe \hqe and this is also true on
a $\mathcal{C}^1$-neighborhood of $J_0$: for any compatible \pce
\ste $J$ satisfying  $|\dd \dd^c_J \psi|<C$ for a constant $C>0$\,
$(M, \omega,J)$ is not \aK-\hq.
\end{thm}
Thus in this case the non-complex \hcte of $(M,J_0)$ is table by
small deformations of the \pce \st. This results is similar to the
proposition \ref{}, but the assumptions on $J$ is much weaker
(\emph{i.e} the open neighborhood is much bigger).\\

The definitions of the notions of "current", "minimizing current" and
"isoperimetric dimension" involved in these theorems will be recalled
in the section \ref{isopsection}. Moreover, in that section we will explain
how the well-bounded geometries allow to cut the \phe curves as wished.
For this, we have proved several isoperimetric results that are presented
in this section. They notably imply:
\begin{lemma}
\label{implicationlemma}
\begin{itemize}
\item Assumption \ref{six} (\emph{resp.} \ref{five}) of theorem \ref{thebigone2}
implies assumption \ref{one} (\emph{resp.} \ref{three}) of theorem
\ref{thebigone1}.
\item Assumption \ref{four} implies assumption \ref{two} in theorem
\ref{thebigone3}.
\end{itemize}
\end{lemma}
Thus in order to prove the three theorems above, one just need to prove
theorem \ref{thebigone1} (then theorem \ref{thebigone2} is just a
corollary), and to prove that theorem \ref{thebigone3} is implied by
its assumption \ref{two}.\\
In fact, we will notice though th proof that these hypothesis of
well-bounded Hessian or radial curvature are needed only to ensure that
the projection on the sets $M_a$ or $B(a)$ are contracting enough. Thus
these assumptions could be replace by ones on the contracting properties of
the projections.\\
\par
These theorems have many consequences in the study of (almost)-complex
hyperbolicities. Let's now state some of them.

\subsection{Applications and examples}
\label{applicexsection}
First of all, let's notice that theorem \ref{thebigone2} obviously implies
Bangert's theorem \ref{bangertthm} on the torus and $\R^{2n}$. Thus it
generalizes
it. Notably, it provids us with a generalization of this result to some
more general manifolds of the form $W=M \exc H$.\\
Keeping the notations of section \ref{applicationw}, the
hypersurface $H$ being
defined by $H=\{s=0\}$, and denoting $\phi=(-log |s|^2)$, $\omega_0=\dd \dd^c
\phi$, $\psi=exp(\phi)$ and $\omega=\dd \dd^c \psi$, then two Stein
manifolds
$(W, \omega_0 , J_0, \phi)$ (non complete) and $(W,\omega, J_0, \psi)$
are at our disposal.\\
Applying theorem \ref{thebigone2} to $(W, \omega, J_0, \psi, |.|)$, one
gets:
\begin{cor}
\label{basicone} Let $(W, \omega, J_0)$ defined as above with $s$ a
pseudo-transverse section. If $(M,\omega)$ is not strongly \spe \hqe
(for example if it is a sub-critical Stein manifold),
and if his Hesssian is well-bounded or his radial curvature is well-bounded,
then (if $m=2$ is an isoperimetric value) for any \pce \ste uniformly
compatible with $\omega$, $(W,\omega,J)$
is not \aKe \hq.\\
Thus, if  $(W. \omega,J_0)$ has a compact quotient $(\td{W},
\td{\omega})$, then $(\td{W},J)$ and $(W,)$ are not Brody-\hqe
for any \pce \ste $J$ compatible with $\td{\omega}$ on $\td{W}$.
\end{cor}
This result can be complemented by considering at the same time the
manifold $(W,\omega_0)$:
\begin{cor}
\label{abitmore}
Let's consider $W$ defined as above with $s$ a pseudo-transverse section, and
let's suppose that $(W, \omega)$ is not strongly \sp-\hqe
\begin{enumerate}
\item if $(W,\omega,J)$ has a well-bounded radial curvature, then for any
\pce \ste $J$ uniformly compatible with  $\omega$ and $\omega_0$,
$(W,J)$ is not \aK-\hqe  (and it is not hyperbolically embedded in $M$),
\item  if $(W,\omega, J, \psi)$ has a well-bounded Hessian, then for any
\pce \ste $J$ uniformly compatible with $\omega$ and $\omega_0$ and
satisfying $|\dd \dd^c_J \phi|_0 \leq A$, with $A$
a constant (and $|.|_0$ the norm associated to $\omega_0(.,J_0.)$),
$(W,J)$ is not Kobayashi-\hq.
\end{enumerate}
Moreover if $J$ can be extended in an \pce \ste on $M$ then either $W$,
either $H$ is not Brody-\hq.
\end{cor}
Let's notice than the condition  $|\dd \dd^c_J \phi|_0 \leq A$
is satisfied by any $J$ close enough to $J_0$ since
$|\dd \dd^c J_0|=1$.\\
A study of the Hessian for some particular manifolds (see \cite{alb}) allow
us to deduce from these corollaries examples where the non-complex \hcte
is stable by small deformations of the \pce \st. Notably:

\begin{cor}
\label{hyperplans}
Let $W= \C \PP^ n \exc \{ k \textrm { hyperplans}\}$, with $k \leq n$,
be provided with its standard \spe and complex \ste $J_0$. Then for any
\pce \ste compatible $J$ in a $\mathcal{C}^1$-neighborhood of $J_0$,
$(W,J)$ is not complex \hq.
\end{cor}
I've got the same result for other manifolds of the form $\PP^n \C \exc
\Sigma$ with $\Sigma$ an algebraic hypersurface satisfying some
technical assumptions. See \cite{alb} for more details. These results
and some more applications will be the issue of a paper to come.\\
\par
Thanks to the result of section \ref{hypproductpart} (about the capacity
of a product being less than the minimum of the two capacities),
the product manifolds constitute
some more examples for which the results of this section can be applied.

\begin{thm}
\label{prodstab}
Let $(M, \omega_M, J_M, \psi_M)$  be a Stein manifold with
$II_\delta$-bounded geometry satisfying some non-\sp-\hcte assumption:
$\mu(M_a, \omega_M) = O\left( a^{(1-\delta)\f{m_M}{2}} \right)$.
Then for any \pce \ste uniformly compatible with
$\omega_M$, $(M,J)$ is not \aK-\hq, and this property is stable by product
with manifold whose geometry is  $II_\delta$-bounded:\\
Let $(N, \omega_N, J_N, \psi_N)$ a Stein manifolds with
$II_\delta$-bounded geometry (with the same $\delta$ as $M$) and which
admits an isoperimetric value bigger than the one of $M$: $m_N \geq m_M$.
If either $M$ or $N$ has a nonpositive radial curvature, then for any
\pce \ste $J$ uniformly compatible with  $\omega =
\omega_M \oplus \omega_N$ (for the metric $g=g_M \oplus g_N$),
the \pce manifold $(M\times N,J)$ is not \aK-\hq.\\
Consequently, for any quotients $M_0$ of $(M, \omega_M, J_M)$ and $N_0$  of
$(N, \omega_N, J_N)$, for any \pce \ste $J$ compatible with
$\omega = \omega_M \oplus \omega_N$ on $M_0 \times N_0$,
the \pce manifolds  $(M_0 \times N_0,J)$
and $(M \times N, J)$ are not Brody-\hq.
\end{thm}
This can be applied with $M=\C^n$, then this theorem reads:
\begin{cor}
Let $(N, \omega_N, J_N, \psi_N)$ be a Stein manifold with
$II_\delta$-bounded geometry and which admits an isoperimetric value
$m_N \geq 2$. Then the product manifold
$\C^n \times N$ is not \aK-\hqe for any \pce \ste $J$ uniformly
compatible with the \spe \ste $\omega=  \omega_0
\otimes \omega_N$ (for the Riemannian metric $g=g_0 \oplus g_N$).\\
Thus for any compact quotient $N_0$ of $(N,\omega_N,J_N)$, for any
\pce \ste $J$ compatible with $\omega_0 \otimes \omega_N$ on
$\T^{2n} \times N_0$, the \pce manifolds $(\T^{2n} \times N_0,J)$
and $(\C^n \times N, J)$ are not Brody-\hq.
\end{cor}
In particular (but the stated result is much more
general), this corollary implies:
\begin{cor}
Let $N$ be a compact hyperbolic manifold, quotient of the \hqe ball
$(\D^n, \omega_{{\rm hyp}})$ (also denoted by $\Hh^n$) by a cocompact
group of \hle isometries. Then for any \pce \ste $J$ compatible with
the product \spe \ste $\omega_0 \oplus\omega_{{\rm hyp}}$ on
$\T^{2n} \times N$, the manifold $(\T^{2n} \times N,J)$ is not Brody-\hq:
there exists a non-constant $J$-\hle map $f:\C \ra \T^{2n} \times N$.
\end{cor}
Indeed, a \hle isometry preserves both the metric and the \pce \st, and
so also preserves the \spe \st. Moreover, the curvature of the
\hqe balls $\D^n$ is strictly negative (a fortiori $K(s)=0$),
their isoperimetric dimension $m=\infty$ (a fortiori they
have isoperimetric values $\geq 2$) and the function $\beta(s) \ra 0$
for $s \ra \infty$ (a fortiori is bounded). So they
satisfy the required assumptions.\\
More generally, for any $M$ satisfying the hypothesis of the theorem
\ref{thebigone2}, not only $(M,J)$ is not \aK-\hqe for any
\pce \ste uniformly compatible, but moreover, this is also true for the
product manifold $M \times \D^n$. Indeed the theorem  \ref{prodstab}
applies with $N= \D^n$ the \hqe ball. Thus,

\begin{cor}
Let $M$ a Stein manifold with $II_\delta$-bounded geometry satisfying the
non-\sp-\hcte assumption: $\mu(M_a, \omega_M) = O\left(
a^{(1-\delta)\f{m_M}{2}} \right)$ with $m$ an
isoperimetric value of $M$. And let $M_0$ be a compact quotient.
Then for any compatible \pce \ste $J$, $M_0$ is not Brody-\hq. Moreover,
for any compact \hqe manifold $H_0$, quotient of the \hqe balls
$\D^n$ by a cocompact group of \hle isometries (for example any
surface of genus $\geq 2$), the product manifold  $M_0 \times H_0$ is not
Brody-\hqe for any \pce \ste compatible with the product \spe \st.
\end{cor}

\par
Before proving the theorems \ref{thebigone1}, \ref{thebigone2}
and \ref{thebigone3}, let's recall the notions of isoperimetric
geometry involved and prove the lemma \ref{implicationlemma}.

\subsection{Isoperimetric notions and lemmas}
\label{isopsection}
Let's recall the definition of $2$-dimensional isoperimetric value of a
riemannian manifold $(M,g)$ following Gromov \cite{gromov3}:
$m$ is a {\bf ($2$-)dimensional isoperimetric value} of $M$ if there exists
a constant $C$ such that, for any contractile loop $l$ in $M$, there exists
a surface $S$ such that $\pp S =l$ and $\area(S)^{1-\f{1}{m}}\leq C \,
{\rm length}(l)$. Besides, the isoperimetric dimension of $M$ is defined
as the upper bound of these isoperimetric values.\\
It is proved in \cite{gromov3} that if $m>2$ then the
manifold is Brody-\hq. That's why in our study of non-\hct, we will
mainly consider some isoperimetric values $m \leq 2$. \\
Let's notice that, according to the definition, if $m$ is an
isoperimetric value of $M \times N$ then it is also an
isoperimetric value for $M$ and $N$. So the isoperimetric dimension
of $M \times N$ is less than the minimum of the ones of $M$ and $N$. And
in fact, one can check that there is equality.

\par
Let's now introduce the other notion occurring in the theorems: the one
of rectifiable $2$-currents. A {\bf{rectifiable $2$-current}} (or more
generally $k$-current) $S$ is a $\R$-linear map on the vector space
of $2$-forms (or more generally $k$-forms) with compact support on $M$.
This notion generalizes the notion of surface in $M$. Indeed, a surface
can be naturally considered as a current which associates to a $2$-form
$\pi$, the number $\int_{S}\pi$.\\
One could find many more details about currents and their properties in
\cite{Federer} (or \cite{Morgan}). Here, in order not to introduce too
many new notations,
we will often replace the vocabulary for currents by the well-known
equivalent vocabulary for the particular case of surfaces. Moreover,
all our considered currents will be $2$-dimensional. So the dimension
will not be specified anymore in the following.\\
The currents which occur in our study are $Q$-minimizing:

\begin{defi}
Let $(M,)$ be a Riemannian manifold and $Q \geq 1$. A current $S$ of $M$
is {\bf $Q$-minimizing} (or {\bf{quasi-minimizing}} with constant $Q$) if
for any Borel set $B$ and any rectifiable current $X$  such that
$\partial(X)=\partial(S \cap B)$,
$\area(S \cap B) \leq Q \area(X)$.
\end{defi}
In the particular case with $Q=1$, this corresponds of the notion of
minimal surfaces. \\
The $Q$-minimizing currents occur naturally in the study of \phe curves,
since, for compatible metrics, the \phe surfaces are $Q$-minimizing:

\begin{lemma}
\label{lem1min}
Let $(M,g)$ be a Riemannian manifold, $J$ be an \pce \ste and $\omega$ be
an exact \spe \ste on $M$. Let's assume that:
\begin{description}
\item[\quad $(i)$] $J$ is uniformly compatible with $\omega$: there
exist some constants $\alpha^2>0$ and $C>0$ such that $\omega(v,Jv) \geq
\alpha^2|v|^2 \; \forall \: v$ and $|Jv| \leq C|v| \; \forall \:v$
\item[\quad $(ii)$] $\omega$ is $g$-bounded: there exists a constant
$\beta>0$ such that $\omega(u,v) \leq \beta^2 |u|\, |v| \;
\forall \: u,v$.
\end{description}
If $f: (S,i) \longrightarrow (M,J)$ is a \phe curve, with $S$ a
compact with boundary Riemannian surface, then the current
$f_{*}(S)$ is $Q$-minimizing with $Q=\frac{\beta^2}{\alpha^2} C.$
\end{lemma}
This is proved for example in \cite{alb} or in \cite{bangert}. Let's notice
that if we pick the \aKe metric $g=\omega(.,J.)$, this result reads as
the well-known result: any $J$-\hle curve is minimal for this metric $g$.\\
This whole section of results on $Q$-minimizing currents has been inspired
by the well-known result for the currents of $\R^{2n}$ (or any $\R^m$):
\begin{prop}
For any constants $Q>1$, $m_0>0$, there exists a constant $c=c(Q,m_0)$
such that: for any $\rho>0$ and any $Q$-minimizing $2$-current in
$\R^{2n}$ with $\partial S \: \subset \: B(c\rho)$ and
$\mathbf{M}(S) \leq m_0 \rho^2$, then
$\mathrm{spt}\,(S) \: \subset \: B(\rho)$.
\end{prop}
This result was used by Bangert in the proof of his theorem \cite{bangert}
in order to cut some \phe curves. In our study, we would need some
equivalent results in a much more general context.\\
First, using the function $\rho$ on the conical ends, we prove:

\begin{prop}
\label{propminriem}
If the radial curvature of a manifold with conical ends $M$ is
well-bounded (and if $m$ is an
isoperimetric value of $M$), then for any constants
$Q>1$ and $m_0>0$, there exists a constant $C=C(Q, m_0)$ such that: for any
$R_0>1$ and for any $Q$-minimizing surface
 $S$ with $\pp S \subset B(C R_0)$ and
$\area (S) \leq m_0 R_0^m$, then $S \subset B(R_0)$ (with $B(r)$ denoting
the set $\{ \rho \leq r \}$).
\end{prop}
In the context of Stein manifold this reads:

\begin{cor}
\label{corminriem}
Let $(M, \omega,J ,\psi)$ be a Stein manifold with
$II_\delta$-bounded geometry, satisfying some non-\sp-\hqe assumption:
$\mu(M_{a}) \leq c \, a^{(1-\delta)\f{m}{2}}$ with $m$ an
isoperimetric value of $M$. \\
Then for any $Q>1$, there exists a constant $C=C(Q)$ such that: for any
$Q$-minimizing surface $S$  with $\pp S \subset
B(C \mu(M_{a+1})^{\f{1}{2}})$ and
$\area (S) \leq \mu(M_{a+1})$ then $S \subset M_a$.
\end{cor}

By following an equivalent approach as the one in proposition
\ref{propminriem}, but in the context of Stein manifold, we could also prove:

\begin{prop}
\label{minimstein}
Let $(M,\omega,J_0, \psi)$ be a Stein manifold with $I_\delta$-bounded
geometry, and let $m$ be an isoperimetric value.
Then, for any $Q>1$, for any $m_0>0$, there exists a constant
$c_0=C(Q,m_0)>0$ such that: for any $Q$-minimizing current $S$ with
 $\pp S \subset M_{c_0 \, a}$ and $\area(S) \leq
m_0 \, {a}^{(1+\delta)\f{m}{2}}$, then $S \subset M_a$.
\end{prop}

The main idea behind these propositions is that, thanks to the assumptions,
the projection on the balls (or respectively the domains $M_a$) are
contracting enough (this explains our remark about these assumptions made
section \ref{enouncethebig}. Let's first detail the proof of proposition
\ref{minimstein}. The other one is then an easiest application of the same
reasoning.\\
Proposition \ref{minimstein} and corollary \ref{corminriem} implies the
lemma \ref{implicationlemma}.\\
\proof[of proposition \ref{minimstein}] Let's consider the projection $p_a$
on $M_a$ along the flow trajectories of $\nabla \psi$. More precisely,
let's $\phi_t$ denote the flow of $X=\f{\nabla \psi}{|\nabla \psi|^2}$,
then for any $x \in M$ and $t \in \R$,
$\psi(\phi_t(x))=\psi(x)+t$. The projection $p_a$ can thus be defined as:
\begin{equation*}
p_a(x)=\begin{cases} x & \textrm{ if } \psi(x)\leq a\\
\phi_{a-C}(x) & \textrm{ if } \psi(x)=C >a
\end{cases}
\end{equation*}
It satisfies:
\begin{lemma}
\label{proj1}
If the Hessian of $\psi$ is well-bounded then for $b>a$,
on $\{ \psi \geq b\}$,
\begin{equation}|(p_a)_*|^2 \leq
  \left({\f{a}{b}}\right)^{2A}
\end{equation}
\end{lemma}

\proof[of lemma \ref{proj1}]
On $\{\psi> a\}$, obviously $(p_a)_*(X)=0$ and
$(p_a)_*(V)= (\phi_{a-b})_*(V)$ for any $V \in T\{\psi=b\}$ (with $b>a$).
Let's notice that if $V \in \{ \dd
\psi=0 \}$, then for any $t$, $(\phi_t)_* V \in \{ \dd
\psi=0 \}$.
Thus, it is to be checked that for any  $V \in T\{\psi=b\}$,
$|(\phi_{a-b})_*(V)|^2 \leq \left({\f{a}{b}}\right)^{2A} |V|^2$.
To this end, let' find a lower bound for $|(\phi_t)_* V|^2$ for $t>0$ and
$V \in T\{\psi= {\textrm{ constant}}\}= \{ \dd \psi=0 \}$. \\
Denoting $\theta(t)=|(\phi_t)_* V|^2$, through a straightforward
computation we check that  $\theta'(t)= \f{2}{|\nabla \psi|^2}Hess
\,\psi(V,V)$. So, $\theta'(t) \geq
\f{2A}{\psi} |V|^2=2A \f{\theta(t)}{t + \psi(x_0)}$ and the integration
of this inequality proves the lemma.
\endproof  \\

Let's now consider a $Q$-minimizing current such that $\pp S
\subset M_{a_0}$ and $\area(S) \leq m_0 \, {a}^{(1+\delta)\f{m}{2}}$.\\
Let's pick any $\gamma >1$ (it will be fixed later) and let's denote
$a_j={\gamma}^{\f{j}{2A}} a_0$, so that on
$\{\psi \geq a_{j+1}\}$, $|(p_{a_j})_*|  \leq \gamma^{-1}$. \\
Then, the currents $S_j=S\cap \{\psi \geq a_j\}$ are $Q$-minimizing and
so satisfy:
$\area(S_j)\leq Q \area(p_{a_j}(S_j))$. Since $|(p_{a_j})_*|  \leq
\gamma$ on $S_{j+1}$, this implies:
\begin{equation*}
(1-\gamma^{-2})\area(S_{j+1}) \leq (1-Q^{-1})\area(S_j).
\end{equation*}
and so at the end, for any $j$,
\begin{equation}
\label{encoreunintermed}
\area(S_{j}) \leq \left(\f{1-Q^{-1}}{1-\gamma^{-2}}\right)^j \area(S_0).
\end{equation}
Moreover, on one side, $\area(S_0) \leq \area(S) \leq m_0
{a}^{(1+\delta)\f{m}{2}}$.\\
And on the other side, let's get a lower bound for $\area(S_{j})$.
Let's suppose that there exists a point $x \in S$ such that
$\psi(x)=R^2$ with $R^2>a_j$. Then we prove:

\begin{lemma}
\label{dimisom}
Let $m$ be an isoperimetric value of $M$. Then there exists a constant
 $C=C(Q,M)$ such that:
for any $Q$-minimizing $S$ in $M$ with $z \in S$ such that  $\{
\psi(z)=R^2 \}$ and $\pp S \subset \{ \psi \leq a \}$, then
\begin{equation*}
\area(S) \geq C\left(\int_{\sqrt{a}}^R \f{1}{\beta_0(s)} ds\right)^m.
\end{equation*}
\end{lemma}

\proof[of the lemma \ref{dimisom}]
Let's consider $a(t)=\area(S \cap \{R-t \leq \sqrt{\psi} \leq R
\})$ and $l(t)=\length(S \cap \{\psi=R-t\})$, then  La coaire formula
writes:
\begin{equation*}
a'(t)\geq \f{2}{\beta_0(R-t)} l(t).
\end{equation*}
Since $S$ is $Q$-minimizing, and by definition of the
isoperimetric value, there exists a constant $c$, depending only on $M$
and $Q$, such that
 $l(t) \geq c \, (a(t))^{1-\f{1}{m}}$. Thus
\begin{equation*}
a'(t)\geq c \, \f{2}{\beta_0(R-t)}\, (a(t))^{1-\f{1}{m}},
\end{equation*}
and the integration proves the lemma. \endproof\\
\indent

According to this lemma \ref{dimisom}, the equation
(\ref{encoreunintermed}) writes:

\begin{equation}
\label{petittruc}
\int_{\sqrt{a_j}}^R \f{1}{\beta_0(s)} ds \leq C
\left(\f{1-Q^{-1}}{1-\gamma^{-2}}\right)^{\f{j}{m}}
(m_0 \, {a}^{(1+\delta)\f{m}{2}})^{\f{1}{m}}
\end{equation}

Because of the assumption of bounded differential, there exist some constants
$c$ and $\delta \geq 0$ such that $\beta(s)\leq c s^\delta$. Since,
$\beta(s)=\beta_0(\sqrt{s})^2$, $\beta_0(s) \leq c' s^\delta$ for some
constant $c'$. So, $\int_{\sqrt{a_j}}^R
\f{1}{\beta_0(s)} ds \geq C'(R^{1+\delta}-(\sqrt{a_j})^{1+\delta})$.
And at the end:
\begin{equation}
\label{petitbis}
R^{1+\delta}-\left( \gamma^{\f{j}{4A}} \sqrt{a_0}\right)^{1+\delta}\leq C
\left(\f{1-Q^{-1}}{1-\gamma^{-2}}\right)^{\f{j}{m}}
m_0^{\f{1}{m}} \, {a}^{(1+\delta)\f{1}{2}}
\end{equation}
In order to conclude the proof of proposition \ref{minimstein}, it remains to
check that there exists a constant $C>0$, independent on $a$, $S$ and $Q$,
such that if $a_0 \leq C\, a$ then (\ref{petitbis})
implies $R \leq a$.\\
To this end let's choose  $\gamma >1$ such that
$\left(\f{1-Q^{-1}}{1-\gamma^{-2}}\right) =\,1-\f{1}{2Q}.$ then,
$\left(\f{1-Q^{-1}}{1-\gamma^{-2}}\right)^{\f{j}{m}}\ra 0$ when $j \ra
\infty$.\\
A $j$ (independent on $a$) can then be fixed so that:
$\left(\f{1-Q^{-1}}{1-\gamma^{-2}}\right)^{\f{j}{m}}
m_0^{\f{1}{m}} \leq 1$. Then,
$R^{1+\delta} \leq (\sqrt{a_j})^{1+\delta}+\f{(\sqrt{a})^{1+\delta}}{2}$
with $a_j =\gamma^{\f{A\, j}{2}}\, a_0$. So with
$c_0 =\left(\f{1}{2}\right)^{\f{1}{1+\delta}} \gamma^{\f{-A\, j}{2}}$, the
proposition \ref{minimstein} is satisfied. \endproof \\

The proof of proposition \ref{propminriem} is similar but uses an approach
based on the distance-function $\rho$.

\proof[of proposition \ref{propminriem}] Let's consider $p_r$ the radial
projection on $B(r)$. Similarly as lemma \ref{proj1}, we check that if
$Hess \rho \geq \f{\mu}{\rho} H$ with $H=g_0 - d\rho \otimes d\rho$, then
$|(p_r)_*| < \left(\f{r}{R}\right)^\mu$ on $\{\rho >R \}$.\\
However as explained before, this hypothesis is satisfied if the radial
curvature is well-bounded.\\
Let's consider any constant $R_0$ and any $Q$-minimizing surface $S$,
satisfying
$\pp S \subset B(r_0)$ et $\area (S) \leq m_0 R_0^m$.\\
Then the same reasoning as above for proposition \ref{minimstein} can be
followed, by first considering the $Q$-minimizing current
 $S_j=S \cap \{\rho>r_j\}$, with  $r_j=
\gamma^{\f{j}{\mu}} r_0$ ($\gamma$ will be fixed later), and this leads to:
\begin{equation*}
\area(S_{j}) \leq \left(\f{1-Q^{-1}}{1-\gamma^{-2}}\right)^j \area(S_0).
\end{equation*}
Let's suppose that there exists $x \in S$ such that $\rho=R$.
Following the same idea than for lemma \ref{dimisom} (but this case is simpler
since $|d\rho|=1$), one gets:
\begin{lemma}
If $m$ is an isoperimetric value of $M$,for any $Q \geq 1$, there exists a
constant $C>0$ such that: for any $Q$-minimizing current $S$ with
$\pp S \subset \{\rho \leq r \}$ and such that there exists
$x \in S \cap \{ \rho = R\}$ with $R>r$, then,
 $\area(S)\geq c(R-r)^m$.
\end{lemma}
Thus there exists a constant $c$ such that for any $j$ (with $r_j \leq R$),
$\area(S_j)\geq c(R-r_j)^m$. Moreover, according to the assumption,
$\area(S_0)\leq \area(S) \leq m_0 R_0^m$. At the end:
\begin{equation}
\label{provisoire}
R \leq \gamma^{\f{j}{\mu}} r_0 + C
\left(\f{1-Q^{-1}}{1-\gamma^{-2}}\right)^{\f{j}{m}} m_0^{\f{1}{m}} R_0.
\end{equation}
We'd like that this inequality could imply $R \leq R_0$ (for an appropriate
choice of $\gamma$, $r_0$...).\\
Let's first fix $\gamma$ so that $\left(\f{1-Q^{-1}}{1-\gamma^{-2}}\right)
= 1 - \f{1}{2Q}$, then $j$ so that $\left(1 - \f{1}{2Q}\right)^{\f{j}{m}}
m_0^{\f{1}{m}} \leq \f{1}{2}$. Then by fixing $r_0=
\f{\gamma^{\f{-j}{\mu}}}{2} R_0$, the inequality (\ref{provisoire}) implies
$R \leq R_0$. So at the end, the proposition is satisfied with
$C=\f{\gamma^{\f{-j}{\mu}}}{2}$ (independent on $R_0$).
\endproof\\
Then the corollary \ref{corminriem} comes directly from the comparison
between the domains $M_a$ and $B(r)$.
\proof[of corollary \ref{corminriem}]

\be
  B\left(\f{2
  s}{\beta_0(s)}\right) \subset M_{s^2}
\ee
Thus, since $\beta(s) \leq c s^\delta$,
\be
\label{comp2}
B(c'a^{\f{1-\delta}{2}}) \subset
B\left( 2 \sqrt{\f{a}{\beta_0(a)}}\right) \subset M_a.
\ee

\indent
To prove the corollary, it needs to be proven that under the made
assumptions,
$S \subset B(R_0)$ avec $R_0=c' a^{\f{1-\delta}{2}}$.\\
Or the assumption states that there exists a constant $c$ such that
$\mu(M_{a+1}) \leq c \,a^{(1-\delta)\f{m}{2}}= c'' \, R_0^m$.\\
Moreover the proposition \ref{propminriem} provided us with a constant
$C_0$ such that for any $R_0>1$, for any $Q$-minimizing current $S$ with
$\pp S \subset B(C_0 R_0)$ and $\area (S) \leq c'' \, R_0^m$ then
 $S \subset B(R_0)$.\\
Let's denote $C=\f{C_0}{c}$. Thus, if  $S$ satisfies $\pp S \subset B(C
\mu(M_{a+1})^{\f{1}{2}})$ and
$\area (S) \leq \mu(M_{a+1})$, then $\pp S \subset B(C_0 R_0^{\f{m}{2}})
\subset B(C_0 R_0)$ et $\area(S) \leq c \, R_0^m$, and so $S \subset B(R_0)$.
\endproof \\

\subsection{Proof of the theorems \ref{thebigone1}
and \ref{thebigone3}}

Let$(M,\omega,J_0, \psi)$ be a Stein manifold, $M_a=\{\psi \leq
a\}$, and $g_0$ be the Riemannian metric associated with
$\omega_0$ and $J_0$. Let's suppose that $(M,\omega)$ is not
strongly \hqe and so for any $a>0$ $\mu(M_a) <\infty$.\\
Let's now consider another \pce \ste  which is uniformly
compatible with $\omega$, which by definition means that there
exist $C>0$ and $\alpha >0$ such that:
\begin{equation}
\label{unif}
\begin{cases}
(i) & \forall \, v, \ |Jv| \leq C|v| \\
(ii) & \forall \, v, \ \omega_0(v, Jv) \geq \alpha |v|^2,
\end{cases}
\end{equation}
Thus the metrics $g=\omega(.,J.)$ and $g_0$ are equivalent.\\
This \pce \ste does not necessarily respect the contact
hyperplane of the $\{ \psi =a\}$. So in order to be able to apply
the theorem \ref{sympdisk} and to get some \phe disks, we are
going to deform the \pce \ste in some new ones equal to $J_0$ on
the neighborhood of infinity (and so which preserves the contact
hyperplanes). This is possible thanks to the contractibility of
the set of the compatible \pce \st s.\\
Indeed, because of the contractibility of $\mathcal{J}=\{ J,
\textrm{\pce \ste on M, compatible with } \omega\}$, there exists
a function $H:[0,1]\ra \mathcal{J}$ such that for any $x \in M$,
$H_0(x)=J(x)$ and $H_1(x)=J_0(x)$, and such that for any $t$ $H_t$
satisfies the properties $(i)$ and $(ii)$ of (\ref{unif}).\\
Then for any $a >0$, we define the function $\lambda_a:\R \ra
[0,1]$, $\smooth$, taking value $0$ for $x \leq a$ and $1$ for $x
\geq a+1$. We can then consider the \pce \st s $J_a \, \in \,
\mathcal{J}$ defined by $J_a(x)=H_{\lambda_a(\psi(x))}(x)$. It
satisfies:
\begin{equation}
\label{jaunif}
\begin{cases}
(i) & \forall \, v, \ |J_Rv| \leq C|v| \\
(ii) & \forall \, v, \ \omega_0(v, Jv) \geq \alpha |v|^2\\
(iii) & J_a=J \textrm{ on } M_a \\
(iv) & J_a=J_0 \textrm{ outside } M_{a+1}
\end{cases}
\end{equation}
Let's notice that the $J_a$ are all uniformly compatible with
$\omega$, and with the same coefficients as for $J$
(\ref{unif}). Moreover, $M_{a+1}$ is $J_a$-convex. Since $\mu(M_a)<\infty$,
if $x_0 \in M$ is fixed, theorem \ref{sympdisk} provides us
with a $J_a$-\hle disk $f_a:\D \lgra M_{a+1}$ with:
\begin{equation*}
\begin{cases}
(i) & f_a(0)=x_0\\
(ii) & f_a(\partial \D) \subset \pp M_{a+1}\\
(iii) & \area_{g_a}(f_a) \leq \mu(M_{a+1})
\end{cases}
\end{equation*}

Let's notice that, since the $J_a$ are uniformly compatible
with $\omega$ for some constants $C$ and $\alpha$ independent on
$a$ (see \ref{jaunif}), the curves $f_a$ are $Q$-minimizing for
a coefficient $Q$-independent on $a$ ($Q=\f{C}{\alpha}$). \\
\par
In order to get a $J$-\hle curve, we need to cut the curves
$f_a$, restricting them to a topological disk $D_a$ so that
$f_a(D_a) \subset M_a$. To this end, we are going to use the
$Q$-minimality of the $f_a$ and the isoperimetric assumptions
of the theorems. \\
\paragraph{Assumption \ref{one} implies the result of theorem
\ref{thebigone1}.}
The $J_a$-\hl s curves $f_a$ being $Q$-minimizing for a constant
$Q$ independent on $a$, let's consider $C=C_Q>0$ the constant provided
by the assumption \ref{one}. (Possibly by replacing $C$ by a generic value
as close to $C$ as wished), $S_a=f_a^{-1}(\Ss_{C \sqrt{\mu(M_{a+1})}})$
is a compact, co-dimension $1$, sub-manifold of $\D$. Let's then $D_a$ be
the biggest topological disk of $\C$ such that
$\pp D_a \subset S_a$ and $0 \in D_a$.
Then the $J_a$-\hle curves ${f_a}_{|D_a}$ satisfy:
\beq*
\begin{cases}
\area(f_a (D_a) \leq \mu(M_{a+1})\\
f_a(\pp D_a) \subset \Ss_{C \sqrt{\mu(M_{a+1})}}\subset
B(x_0,C \sqrt{\mu(M_{a+1})}),
\end{cases}
\eeq*
and so, according to assumption \ref{one}. $f_a(D_a)\subset M_a$.
Consequently the curves ${f_a}_{|D_a}$ are $J$-\hl.\\
By reparametrizing $f_a$ via a biholomorphism from
$\D$ on $D_a$, we get some $J$-\hle curves $h_a: \D \lgra M_a$ with:
\begin{equation*}
\begin{cases}
(i) & \area(h_a) \leq \mu(M_{a+1}) \\
(ii) & h_a(0)=x_0\\
(iii) & h_a(\partial \D) \subset \Ss(x_0,C \sqrt{\mu(M_{a+1})})
\end{cases}
\end{equation*}
Then, $\area(h_a)=O\left(d^2(h_a(0), h_a(\partial \D)) \right)$ and
the proposition \ref{lemmebgt} implies that $(M,J)$ is not \aKe-\hq.
\paragraph{Assumption \ref{three} implies the result in theorem
\ref{thebigone1}.}
Again, as the $f_a$ are all $Q$-minimizing for the same constant $Q$
independent on $a$, let's consider the constant $C=C_Q$ provided by the
assumption \ref{three}. Then we apply the same usual reasoning: by
first considering the codimension $1$ submanifold
$S_a = f_a^{-1}(\{\psi =C\, a \})$ (or possibly
$f_a^{-1}(r_a)$ for some $r_a$ close to $C \,a$ by inferior values),
then $D_a$ the biggest topological disk such that $\pp D_a \subset S_a$
and $0 \in D_a$, we get some $J$-\hle disks $h_a:\D \ra M_a$
with:
\begin{equation*}
\begin{cases}
(i) & \area(h_a) \leq \mu(M_{a+1}) \\
(ii) & h_a(0)=x_0\\
(iii) & h_a(\partial \D) \subset \{\psi=C \, a \}.
\end{cases}
\end{equation*}
Besides, one can notice that $\f{2s}{\beta_0(s)} \leq
d(x_0,\{\phi=s^2\})\  (+C)$. Indeed if $\sqrt{\psi(x)} =s$ and
$\gamma$ is a path between $x$ and $x_0$ then
\beq*
s \leq \sqrt{\psi(x)}- \sqrt{\psi(x_0)}= \int_\gamma \f{\dd \psi}{2 \sqrt{\psi}}
\leq \f{\beta_0(s)}{2} {\length}(\gamma) \textrm{ and }  B\left(\f{2
  s}{\beta_0(s)}\right) \subset M_{s^2}.
\eeq*
So, $d(x_0,\{\psi=C \, a\}) \geq
2 \sqrt{\f{C \, a}{\beta(C \, a)}}$. Thus if there exists a sequence
$(a_n)$ with $\mu(M_{a_n+1})=O\left(\f{C \, a_n}{\beta(C \, a_n)}\right)$,
the proposition \ref{lemmebgt} can be applied and $(M, \omega,J)$ is
not \aK-\hq.

\paragraph{Assumption \ref{two} implies the result in theorem
\ref{thebigone3}}
The isoperimetric assumption being the same as the one in \ref{three}
of theorem \ref{thebigone1}, one can begin as above and get the same
sequence $(h_a)$ of \phe diks. However then, our assumption of
non-\spe \hcte is less strong and we cannot apply the same reasoning.\\
However if we suppose that $J$ satisfies $|\dd \dd^c_J \psi| <C$ for
a constant $C$, the norm being taken for the metric $g_0$ (since
$|\dd \dd^c_{J_0} \psi|=1$ this will be notably the case if $J$ is
 $\mcl{C}^1$-close enough to $J_0$), then $\tau(h_a)\geq C' a$.
Indeed, with a proof similar to the one of lemma \ref{taupsi}, one
can check the equality ($h_a$ is $J$-\hl):
\begin{equation*}
\int_{0}^{2\pi}\psi \circ h_a(e^{i \theta}) \dd\theta \,-\, 2\pi
\psi\circ h_a(0)=\int_{0}^{1} \frac{d\rho}{\rho} \int_{\D(\rho)} h_a^*
\dd\dd^c_J \psi.
\end{equation*}
Moreover, because of our assumption:

\begin{equation*}
\int_{0}^{1} \frac{d\rho}{\rho} \int_{\D(\rho)} h_a^* \dd \dd^c_J \psi
\leq C \, \tau(h_a),
\end{equation*}
and thus, $\tau(h_a) \geq C' a$.
Consequently, if $\mu(M_{a_n})=O(a_n)$ for an exhaustive sequence $(a_n)$,
then the proposition \ref{proptau} can be applied and $(M,\omega,J)$ is
not \aK-\hq.\\
\par
This concludes the proof of theorem \ref{thebigone3}, since
 assumption \ref{four} implies assumption \ref{two}.
And according to lemma \ref{implicationlemma}, this also proves theorem
\ref{thebigone2}.
\par
Now it just remains to prove the corollaries of section
\ref{applicexsection}.

\subsection{Proof of the theorems of section \ref{applicexsection}}
\paragraph{Proof of corollaries \ref{basicone} and \ref{abitmore}}
First let's consider a manifold of the type  $W=M \exc H$, keeping the usual
notations. Moreover let $|.|$ denote the metric associated to $\omega$
and $J_0$,  and $|.|_0$ the one associated with $\omega_0$ and $J_0$. \\
Let's notice that
\be
\label{difsp}
\omega = \psi \ \omega_0 + \psi \ \dd \phi \wedge \dd^c \phi
\ee
With a straightforward computation, this implies
\begin{equation*}
|\dd \psi|^2 = \psi \, \f{|\dd\phi|_0^2}{1+|\dd\phi|_0^2}.
\end{equation*}
However, according to the inequality (\ref{inegverifnormesigma}),
if the section $s$
defining $H$ is pseudo-transverse, then
$|\dd \phi|_0 \geq c \, \psi^{1- \alpha}$ and,
\be
\label{primsestim}
\f{\psi^{1-\alpha}}{C + \psi^{1-\alpha}} \leq \f{|\dd \psi|^2}{\psi} \leq 1.
\ee
This provides us with some estimates for the fundamental functions
$\alpha$ and $\beta$. First, $\beta(s) \leq 1$ (the assumption of
theorem \ref{thebigone2} is satisfied with $\delta=0$).\\
Moreover, $\f{1}{\alpha(\sqrt{s})^2} \leq 1+ \f{1}{s^{1-\alpha}}$.
Thanks to the estimate of the capacity  in (\ref{dansunsens}), there
exists a constant $C$ such that $\mu(M_a, \omega) \leq C \, a$
(with $M_a=\{ \psi \leq a \}$). Let's emphasize that we've proved
(comparing also with approximation already got for the capacity for
the \spe \ste $\omega_0$):
\begin{lemma}
\label{cappseudotrans2}
Considering $W=M \exc H$ with $H$ defined as the zero set of
a pseudo-transverse
section $s$ of a hermitian \hle line bundle $(\mcl{L},|.|, \nabla)$
over the \Kae manifold $(M, \omega_0, J_0)$ with curvature form $\omega_0$,
and denoting $\phi=(-log |s|^2)$, $\omega_0=\dd \dd^c
\phi$, $\psi=exp(\phi)$ and $\omega=\dd \dd^c \psi$, if $W$ is not \spe
hq, then
\begin{itemize}
\item $\forall a>0$, there exists $C>0$ such that
$\mu(\{\phi \leq a \}, \omega_0) \leq C \ a^\alpha$
\item there exists $C>0$ such that $\mu(\{ \psi \leq a \}, \omega) \leq C a$
\end{itemize}
\end{lemma}
That's why, if the assumptions of corollary  \ref{basicone} are
satisfied then the assumption \ref{five}. or \ref{six}. of
theorem \ref{thebigone2} is satisfied and this proves corollary
\ref{basicone}.\\
\par
In order to prove corollary \ref{abitmore}, the assumptions, (and so also the
conclusions) of \ref{basicone} are kept.\\
Since the assumption  \ref{five}. or \ref{six}. of
theorem \ref{thebigone3} is satisfied by $(W,J,
\omega, \psi)$, one gets, for any $J$ uniformly compatible with $\omega$,
a sequence of $J$-\hle curves (by following the demonstration of
theorem  \ref{thebigone1} or \ref{thebigone2}):
$ f_a: \D \ra M_a$, satisfying either:
\begin{equation*}
\begin{cases}
(i) & \area(f_a) \leq \mu(M_{a+1}) \\
(ii) & f_a(0)=x_0\\
(iii) & f_a(\partial \D) \subset \{\psi=C \, a \},
\end{cases}
\end{equation*}
if the Hessian is well-bounded, either:
\begin{equation*}
\begin{cases}
(i) & \area(f_a) \leq \mu(M_{a+1}) \\
(ii) & f_a(0)=x_0\\
(iii) & f_a(\partial \D) \subset \Ss(C \, \sqrt{a}),
\end{cases}
\end{equation*}
if the radial curvature is well-bounded. \\
Actually, it had just been proven that  $f_a(\partial \D) \subset
\Ss(x_0, c \sqrt{\mu(M_{a+1})})$, but in this particular case,
we've checked that  $\mu(M_a) \leq C \, a$ (with $C$ constant).

\par
Moreover, if $J$ is supposed to be uniformly compatible with
$\omega_0$, then, using Stokes' theorem, one gets
\begin{equation}
\area(f_a, g_0)\leq C \int_\D f_a^* \omega_0 = C' \f{1}{a} \int_{\D}
f_a^* \omega \leq C'' \f{1}{a} \mu(M_a) \leq C_0.
\end{equation}
In the case of the radial curvature assumption, this provides us directly
with $\area(f_a, g_0)= o(d_{g_0}(f_a(\pp \D),f_a(0)))$. So according to
proposition \ref{lemmebgt}, $(W,J, g_0)$ is not \aKe \hqe (and more precisely,
$g_0$-\hq). Since $(W,g_0)$ is relatively compact into
$M$, this implies that  $(W,J)$ is not hyperbolically embedded into
$M$. The last part of the corollary is then obvious.
\par
In the case of the Hessian assumption, if we suppose moreover that
 $|dd \dd^c_J \phi|_0 \leq A$ (with $A$ any constant), it has already been
 proven in the last section that:
\begin{equation*}
\tau(f_a,\omega_0) \geq B \inf_{f_a(\pp \D)} \phi \circ f_a,
\end{equation*}
and so in this case $\tau(f_a) \geq C' \, a$. Then,
$\area(f_a, g_0)= o(\tau(f_a,\omega_0))$. And proposition
\ref{proptau2} provides us with the wished conclusion.\\

\paragraph{Proof of theorem \ref{prodstab}}
First let's notice that if $m_M=2$ is an isoperimetric value if $M$ and if
$N$ has one $m_N \geq 2$, then $2$ is an isoperimetric value of
$M \times N$.\\
Moreover by definition of the radial curvature one can easily check that
for any $(x,y) \in M \times N$, the radial curvature of the product
$M \times N$ satisfies :
\beq*
curv_{M \times N} (x,y) \leq \max(curv_M(x), curv_N(y)).
\eeq*
Indeed, the metric considered being the product-metric,
the radial vector of  $M \times N$,
 $\pp_{M\times N} =\f{1}{\sqrt{2}}(\pp_M +\pp_N)$ and the curvature tensor
 is the sum of the curvatures tensors of $M$ and $N$.
Thus, $ \dsps curv_{M \times N} (x,y) =
\f{1}{\sqrt{2}} \sup_{\overset{v_1 \in TM, v_2 \in TN}{|v_1|^2+|v_2|^2=1}}
(|v_1| curv_M(x)+|v_2| curv_N(y))$, and the Cauchy-Shwarz inequality proves
the inequality. In particular, if $K_M=0$ then :
\beq*
K_{M \times N}(s) \leq K_N(s).
\eeq*
Thus if one of the radial curvature is null and the other one satisfies the
inequality $\int s \, K(s) <1$, then the one of $M \times N$ satisfies this
same inequality. The proof is then similar to the one of the theorems
\ref{thebigone1} or \ref{thebigone2}. Let's consider the $J_0$-convex sets
$M_a=B_M(a) \times B_N(a)$ (with $J_0=J_M \times J_N$) and their
capacity $\mu(M_a, \omega) \leq \mu(B_M(a),
\omega_M) \leq C a^2$ (with $\omega=\dd \dd^c \rho_M^2 \otimes \omega_N$). \\
For any fixed uniformly compatible $J$ with $\omega$  fix{\'e}, we build
some $J_a$ by deforming $J$ in $J_0$ at infinity
($=J$ on $M_a$, $J_0$ outside $M_{a+1}$) and we get some $J_a$-\hle maps
with $\area(f_a) \leq C a^2$ and $f_a(\pp \d) \subset \pp M_{a+1}$, that
need to be cut. For this we use proposition \ref{propminriem}, that have
been stated for the sets $B_{M\times N}(R)$. But since
 $B_M(\f{a}{\sqrt{2}}) \times B_N(\f{a}{\sqrt{2}})
\subset B_{M \times N} (a) \subset B_M(a) \times B_N(a) =M_a$, we get as a
corollary an identical result for the $M_R$. And finally, this provides us
with sequence of $J$-\hl-maps $h_a :\D \ra M \times N$ with
$\area(h_a) \leq C a^2$ and
$h_A(\pp D) \subset \pp M_{c_0 \, a}$. Proposition \ref{lemmebgt} concludes
the proof of theorem \ref{prodstab}.

\end{document}